\documentstyle[12pt]{article}

\newtheorem{theorem}{Theorem}[section]

\newtheorem{lemma}{Lemma}[section]
\newtheorem{corollary}{Corollary}[section]
\newtheorem{prop}{Proposition}[section]
\newtheorem{definition}{Definition}[section]
\newtheorem{remark}{Remark}[section]
\newtheorem{example}{Example}[section]
\begin{document}
\title{On some generalizations of\\ Jacobi's Residue Formula.
\footnote {AMS classification 
number: 32A27, (32A25, 32C30).}}
\author { by \\ Alekos Vidras \\ Department of Mathematics and Statistics\\
           University of Cyprus, Nicosia, Cyprus \\
         Alain Yger, \\
         Department of Mathematics, \\
         University of Bordeaux, Talence,  France.  }
 \maketitle
 \begin {abstract}
Using Bochner-Martinelli type residual currents 
we prove some generalizations of 
Jacobi's Residue Formula, which allow proper polynomial 
maps to have  'common zeroes at  infinity  ', in projective or toric 
situations. 
\end{abstract}
\section{Introduction}
\setcounter{equation}{0}
One of the classical results in the one complex variable 
residue theory is the following: 
for every polynomial map $P:{\bf C}\rightarrow {\bf C}$, 
the total sum of residues of the form ${Qd\zeta/P}$ (where 
$Q\in {\bf C}[X]$) at the zeroes of
$P$ equals the residue at infinity of the rational 
function $P/Q$ with the opposite sign. \\
Some multidimensional analogues of this statement are treated 
in the present note. Consider the polynomial 
map
\begin{eqnarray*}
P=(P_1,\dots,  P_n) :{\bf C}^n\longrightarrow {\bf C}^n
\end{eqnarray*}
and assume that ${\bf C}^n$ is imbedded into the 
complex projective space ${\bf P}^n$. 
Let $^hP_1, \dots , ^hP_n $ be the homogenizations 
of the $P_j$, $j=1,\dots , n$, that is the homogeneous polynomials in $n+1$ 
variables 
\begin{eqnarray*}
^hP_j(X_0, X_1, \dots , X_n)=X_0^{\deg P_j}
P_j({X_1\over X_0}, \dots , {X_n\over X_0}).
\end{eqnarray*}
Let us impose the Jacobi condition, that is 
\begin{eqnarray} 
& &\hbox{{\em The homogeneous parts of higher degree in }}\; 
P_j(X_1, \dots , X_n), \; {\em for} \nonumber \\
& &j=1,\dots,n,\ \hbox{{\em do not have common zeroes in}}
 \;{\bf C}^n\setminus (0,\dots,0). 
\end{eqnarray} 
Then, it is a clasical result that 
goes back to Jacobi \cite{j:gnus}, that the set 
$V(P):=\{P_1=\dots=P_n=0\}$ is finite, 
with cardinal number equal to $\deg P_1\cdots \deg P_n$ and that 
for any $Q\in {\bf C}[X_1,\dots, X_n]$,
 such that 
 \begin{eqnarray*}
 \deg Q\leq \sum\limits_{j=1}^n \deg(P_j)-n-1,
 \end{eqnarray*}
 one has 
 \begin{eqnarray}
 {\rm Res}\, \left [
 \begin{array}{ccccc}
 Q(X_1,\cdots, X_n)dX\\
 P_1,\cdots ,P_n\\
 \end{array}
 \right ] 
 = \sum\limits_{\alpha \in V(P)}
{\rm Res}_{\alpha }[{Qd\zeta \over P_1\cdots P_n}]=0\,.
 \end{eqnarray}
Here $dX$ means as usual for $dX_1\wedge \cdots \wedge dX_n$ and  
the residue of the meromorphic form ${Qd\zeta \over P_1\dots P_n}$  
at the  isolated point $\alpha\in \{P_1=\dots=P_n=0\}$ is defined as 
 \begin{eqnarray*}
 {\rm Res}_{\alpha }[{Qd\zeta \over P_1\dots P_n}]
={1\over (2\pi i)^n }\lim\limits_{\begin{array}{ccccc}
 \epsilon_1 \mapsto 0\\
 \cdots \\
 \epsilon_n \mapsto 0\\
 \end{array} } \int \limits_{{}_{\begin{array}{ccccc}
 \vert f_1\vert =\epsilon _1\\
 \cdots \\
 \vert f_n \vert =\epsilon_n\\ \zeta \in U_\alpha \\
 \end{array} }} {Q(\zeta )d\zeta \over P_1(\zeta )\dots P_n(\zeta )},
 \end{eqnarray*}
where $U_{\alpha }$ is any bounded domain in ${\bf C}^n$ such that 
$\{\alpha\}=U_\alpha \cap \{P_1=\cdots=P_n=0\}$ and the orientation for the 
cycle $\{\zeta \in U_\alpha,\ |f_1|=\epsilon_1,\dots,|f_n|=\epsilon_n\}$ is 
the one that respects the positivity of the differential form $d\,{\rm arg} 
(f_1)\wedge \cdots \wedge d\, {\rm arg} (f_n)$. 
\vskip 2mm
\noindent
The result of Jacobi has a toric pendant which is 
due to A. Khovanskii \cite{kh:gnus}. 
Let ${\bf T}^n=({\bf C}^*)^n$ and 
$F_1, \dots F_n$ be $n$ Laurent polynomials in $n$ variables
\begin{eqnarray*}
F_j(X_1,\dots X_n)=\sum\limits _{\alpha _j\in {\cal A}_j}
c_{j,\alpha _j}X_1^{\alpha _{j1}}\dots X_n^{\alpha _{jn}},\ j=1,\dots , n, 
\end{eqnarray*} 
with $c_{j,\alpha_j}\not=0$ for 
any $j\in \{1,\dots,n\}$, any $\alpha_j\in {\cal A}_j$ 
(the ${\cal A}_j$ are the supports of the $F_j$). Let 
$\Delta_j$ be the Newton polyhedron of $F_j$, which is 
by definition the closed 
convex hull of ${\cal A}_j$, $j=1,\dots,n$. We now impose the 
Bernstein condition \cite{be:gnus}, that is
\begin{eqnarray} 
& &\hbox{{\em For any }}\  \xi \in {\bf R}^n \setminus (0,\dots,0), 
\hbox{{\em the intersection with}}\ {\bf T}^n \ 
\hbox {{\em of the set}}\; 
\nonumber \\
& &\Big\{\zeta;\ \sum\limits_{{{
 \alpha _j\in {\cal A}_j}\atop {
 <\alpha _j, \xi >=\min\limits_{\eta \in \Delta _j }<\eta, \xi >}}}
  c_{j,\alpha _j}\zeta_1^{\alpha _{j1}}\dots \zeta_n^{\alpha _{jn}}=0, 
j=1,\dots , n \Big\}\ \hbox{{\em is empty}}.  
\end{eqnarray} 
Under such hypothesis, D. Bernstein proved in 
\cite{be:gnus} that the set 
$V^*(F):=
\{F_1=\dots =F_n=0\}\cap {\bf T}^n$ is finite with cardinality 
equal to the Minkowski mixed volume of $\Delta_1,\dots , \Delta_n$ and 
A. Khovanskii \cite{kh:gnus} proved that for any Laurent polynomial $Q$ 
whose support lies in the relative interior of the convex polyhedron 
$\Delta_1+\cdots+\Delta_n$ (that is the interior of this polyhedron 
in the smallest affine subspace of ${\bf R}^n$ that contains it), one has
\begin{eqnarray}
{\rm Res}\left [
 \begin{array}{ccccc}
 Q(X_1,\cdots, X_n)dX\\
 F_1,\cdots ,F_n\\
 \end{array}
 \right ]_{{\bf T}} 
 :=\sum\limits_{\alpha \in V^*(F)}
{\rm Res}_{\alpha }[{Q \over F_1\cdots F_n}
 {d\zeta\over \zeta _1\dots \zeta _n}]=0\,.
\end{eqnarray}
\vskip 3mm
We will see in section 2 how 
it is essential to interpret both geometrically  and analytically 
the conditions (1.1) imposed on $(P_1, \dots , P_n )$ in 
the projective setting or the 
conditions (1.3) imposed on $(F_1, \dots , F_n)$ in the toric setting.\\

In the first case (that is the projective one), the set of conditions (1.1)  
is geometrically  equivalent to the fact that the $n$ 
Cartier divisors ${\cal D}_1, \dots, {\cal D}_n$, defined on ${\bf P}^n$ 
by the homogeneous polynomials 
$^h P_j(X_0,\dots , X_n)$, $j=1,\dots,n$, are such 
that their supports $|{\cal D}_j|$ satisfy 
\begin{eqnarray*}
|{\cal D}_1|\cap \dots \cap |{\cal D}_n|\subset {\bf C}^n. 
\end{eqnarray*}
From the analytic point of view, this is equivalent to the 
following strong properness condition on the polynomial map 
$P=(P_1,\dots,P_n)$ from ${\bf C}^n$ to ${\bf C}^n$: 
there are constants $R>0$, $c>0$, such that, for $\|\zeta\|\geq R$, 
\begin{eqnarray}
\sum\limits _{j=1}^n{\vert P_j(\zeta) \vert\over (1+\| \zeta \|^2)^
{{\deg P_j\over 2}}}\geq c \,. 
\end{eqnarray}
  
  In the toric case, given a smooth toric variety ${\cal X}$ 
associated to any fan which is a simple refinement of the fan 
attached to the polyhedron $\Delta _1+\dots +\Delta _n$, 
conditions $(1.3)$ mean that the effective Cartier divisors
\begin{eqnarray*}
{\cal D}_j= {\rm div} (F_j)+E(\Delta _j), 
\end{eqnarray*}
where $E(\Delta_j)$ is the $\bf T$-Cartier divisor on 
${\cal X}$ associated with the polyhedron $\Delta _j$ (see 
\cite {fu:gnus}), are such that 
\begin{eqnarray*}
\vert {\cal D}_1\vert \cap \dots \cap \vert {\cal D}_n\vert  
\subset {\bf T}^n.
\end{eqnarray*} 
The analytic interpretation of this is the following: 
there exist constants $R>0, c>0 $ such that, for $\zeta\in 
{\bf C}^n$ such that $\|{\rm Re}\, \zeta\|\geq R$, 
\begin{eqnarray}
\sum\limits_{j=1}^n 
{\vert F_j(e^{\zeta_1}, \dots , e^{\zeta_n})
\vert 
\over  e^{H_{\Delta _j}({\rm Re}\, \zeta )}}\geq c\,,
\end{eqnarray}
where $H_{\Delta _j}$ denotes 
the support function of the convex polyhedron $\Delta _j$, that is 
the function from ${\bf R}^n$ to ${\bf R}$ defined as follows
\begin{eqnarray*}
H_{\Delta _j}(x ):=\sup \limits_
{\xi \in \Delta _j }<x, \xi >, \; x\in {\bf R}^n.
\end{eqnarray*}
\vskip 2mm
In \cite {by1:gnus, by2:gnus, by3:gnus}, 
one used extensively the fact that an analogous version 
of $(1.2)$ could be stated whenever the polynomial map 
\begin{eqnarray*}
P=(P_1,\dots ,P_n): {\bf C}^n\mapsto {\bf C}^n
\end{eqnarray*}
was proper. 
We will prove in section 3 of this paper what appears 
to be the sharpest version of such a result, namely
\begin{theorem}
Let $P=(P_1,\dots ,P_n)$ be a polynomial map 
from ${\bf C}^n $ to ${\bf C}^n$ such that there exist 
constants $c>0, R>0$, and rational numbers 
$0<\delta_j \leq \deg(P_j)$, $j=1,\dots , n$, in order that, for 
$\|\zeta\|\geq R$,  
\begin{eqnarray}
\sum\limits_{j=1}^n{\vert P_j(\zeta)\vert 
\over (1+\|\zeta\|^2)^{{\delta _j\over 2}}}\geq c\,.
\end{eqnarray}
Then, for any polynomial $Q\in {\bf C}[X_1,\dots, X_n]$ 
which satisfies 
\begin{eqnarray*}
\deg Q\leq \delta _1+\dots +\delta _n -n-1,
\end{eqnarray*}
one has 
\begin{eqnarray}
{\rm Res}\,\left [
\begin{array}{ccccc}
Q(X_1,\cdots, X_n)dX\\
P_1,\cdots ,P_n\\
\end{array}
\right ]=0.
\end{eqnarray}
\end{theorem}
We will also prove in the same section the corresponding toric version, namely 
\begin{theorem}
Let $F=(F_1,\dots ,F_n)$ be 
a system of Laurent polynomials in n variables, with 
respective Newton polyhedra $\Delta_1,\dots,\Delta_n$.  
Suppose there exist constants $c>0, R>0$, and convex polyhedra 
$\delta_1,\dots ,\delta_n $ with vertices in ${\bf Q}^n$, with 
$\delta_j\subset  \Delta _j$, $j=1,\dots,n$ and 
$\dim(\delta_1+\dots +\delta _n)=n$, which are such that, for any $\zeta\in 
{\bf C}^n$ with $\|{\rm Re}\, \zeta\|\geq R$,
\begin{eqnarray} 
\sum\limits_{j=1}^n{\vert F_j(e^{\zeta_1},\dots ,e^{\zeta_n})
\vert \over e^{H_{\delta _j}({\rm Re}\, \zeta )}} \geq c\,.
\end{eqnarray}
Then, for any  Laurent polynomial $Q$ such that the support of $Q$ 
lies in the interior of the convex polyhedron $\delta_1+\cdots+\delta_n$, one
has
\begin{eqnarray}
{\rm Res}\,\left [
\begin{array}{ccccc}
Q(X_1,\cdots, X_n)dX\\
F_1,\cdots ,F_n\\
\end{array}
\right ]_{{\bf T}}=0.
\end{eqnarray}
\end{theorem}
 
The main tool to be used in the proofs of both theorems 
will be the Bochner-Martinelli integral formula suitably adapted to each case.
 
From the point of view of algebraic geometry 
such theorems are not classical in nature since the supports of the Cartier 
divisors ${\cal D}_1,\dots, {\cal D}_n$ on ${\bf P}^n$ 
corresponding to the $^hP_j$ 
in the first case, or the supports of the divisors 
${\cal D}_j={\rm div} (F_j)+E(\Delta_j)$ on a convenient smooth 
toric variety ${\cal X}$ in the second case, do not intersect properly 
in ${\bf P}^n$ or in ${\cal X}$ (the intersection is assumed to be 
proper in ${\bf C}^n$ or in ${\bf T}^n$).  
Following the point of view developped by J. Koll\'ar in \cite{ko1:gnus, 
 ko2:gnus}, 
or by Lazarsfeld-Ein in \cite{elaz:gnus}, we will also present in section 3 
a geometric interpretation of the conditions $(1.7)$ (in the projective case) 
or $(1.9)$ (in the toric case). 
We will see that the Bochner-Martinelli representation
formula we used below fits 
with the construction of residue currents 
in the non-complete intersection case, 
as proposed in \cite{pty:gnus}. 
A better understanding of our two theorems will lie 
then on the fact that, if $f_1,\dots,f_n$ are $n$ 
holomorphic functions in some domain $\Omega$ of 
${\bf C}^n$, a crucial property 
of the distribution $T_f\in {\cal D}'(\Omega)$ 
whose action on 
a test form $\varphi \in {\cal D}(\Omega)$ is defined (see 
for example \cite{pty:gnus}) 
by
\begin{eqnarray*}
T_f(\varphi):=\lim\limits_{\epsilon \mapsto 0}\ 
{1\over \epsilon ^n}
\int\limits_{\vert f_1\vert ^2+\dots +\vert f_n\vert^2=\epsilon }
\sum\limits_{k=1}^n
(-1)^{k-1}\bar f_k\bigwedge\limits_{{l=1}\atop{l\not=k}}^n 
\overline {\partial f_l} \wedge \varphi d\zeta,
\end{eqnarray*}
is that it is annihilated, as a distribution, 
by any holomorphic function in $\Omega$ which is locally in the 
integral closure of the ideal $(f_1,\dots,f_n)^n$ (this ideal 
is contained in $(f_1,\dots,f_n)$ by the classical result 
of Brian\c con-Skoda \cite{bs:gnus}). Therefore, once 
the  hypothesis will be settled in a natural geometric context, our 
two theorems will appear to be in close relation with 
this  Brian\c con-Skoda theorem, which also 
plays a significant role in \cite{ko2:gnus}, \cite{elaz:gnus}, as 
a transition tool between Lojasiewicz inequalities (or regular 
separation conditions) and effective versions of the algebraic 
Nullstellensatz. 
    
  As a consequence, it will be then natural to present in section 4 
some applications of our two theorems to effectivity questions related to 
the algebraic Nullstellensatz in the projective case or the sparse 
Nullstellensatz in the toric case,  under some properness 
assumptions on the data in ${\bf C}^n$ or in ${\bf T}^n$. 
Such results will extend or sharpen some previous 
results in \cite{by2:gnus, by3:gnus, fpy:gnus, y:gnus}. 
We will also suggest possible applications to some 
results of Cayley-Bacharach type (see \cite{egh:gnus}), 
in the context of improper intersections on ${\bf P}^n$ 
or on a smooth toric variety ${\cal X}$.

{\bf Acknowledgments.} We are indebted to Institut Culturel Francais in Cyprus and the University of Cyprus for the financial support 
they provided during the preparation of this work.\\

\section{An analytic interpretation of Jacobi 
\hfil\break
or Bernstein conditions}
\setcounter{equation}{0}
Using the notation of the previous section we will state in analytic terms 
the conditions (1.1) or (1.3).
    We begin with the 
    \begin{prop}
 Let $P_1,\dots , P_n $ be $n$ polynomials in ${\bf C}[X_1, X_2, \dots, X_n]$.
 The following two assertions 
 are equivalent\\
i) $\{\zeta\in {\bf C}^{n+1},\ 
{}^h P_1=\dots ={}^h P_n=\zeta_0=0\}=\emptyset $\\
ii) There exist strictly positive constants $R,c $ such that, for any 
$\zeta\in {\bf C}^n$ with $\|\zeta\|\geq R$, 
\begin{eqnarray}
\sum\limits_{j=1}^n{\vert P_j(\zeta)\vert 
\over (1+\|\zeta\|^2)^{{\deg P_j\over 2}}}\geq c .
\end{eqnarray}
\end{prop}
{\bf Proof.} Writing {\sl ii}) in homogeneous coordinates, we get that, 
if $(\zeta_0,\dots,\zeta_n)\in {\bf C}^{n+1}$ is such that 
$$
\vert \zeta_1\vert +\dots +\vert \zeta_n\vert >R \vert \zeta_0\vert\,, 
$$
one has
$$
\sum\limits_{j=1}^n\vert {}^h P_j(\zeta_0, \zeta_1, \dots ,\zeta_n)\vert
 \geq c(\sum\limits_{j=1}^n (\vert \zeta_0\vert ^2+\dots +\vert 
\zeta_n\vert ^2)^{{\deg P_j\over 2}}.
$$
In particular 
$$
\sum\limits_{j=1}^n\vert {}^h 
P_j(0, \zeta_1, \dots ,\zeta_n)\vert
 \geq c(\sum\limits_{j=1}^n (\vert \zeta_1\vert ^2+\dots +\vert \zeta_n
\vert ^2)^{{\deg P_j\over 2}}.
$$
This shows that ii) implies i).\\
Let now $P_j(X)=p_j(X)+q_j(X)$, such that $\deg q_j<\deg p_j $, 
$p_j$ being an homogeneous polynomial with degree 
$d_j=\deg (P_j)$ (the leading terms in $P_j$). 
Condition {\sl i}) is equivalent to the fact that 
\begin{eqnarray*}
\{\zeta \in {\bf C}^n,\ p_1(\zeta)=\dots =p_n(\zeta)=0\}=\{(0,\dots,0)\}.
\end{eqnarray*}
Since $p_1,\dots p_n $ are homogeneous with respective 
degrees $d_1, \dots , d_n $, there exists $c>0$ such that, for any 
$\zeta\in ({\bf C}^n)^*$, 
\begin{eqnarray*}
\sum\limits_{j=1}^n{\vert p_j(\zeta)\vert \over \|\zeta\|^{d_j}}>
\tilde c\,.
\end{eqnarray*}
Therefore, for any $\zeta\in ({\bf C}^n)^*$, one has 
\begin{eqnarray*}
\sum\limits_{j=1}^n{\vert P_j(\zeta)\vert \over \|\zeta\|^{d_j}}\geq 
\sum\limits_{j=1}^n{\vert p_j(\zeta)\vert \over \|\zeta\|^{d_j}}
-\sum\limits_{j=1}^n{\vert q_j(\zeta)\vert \over \|\zeta\|^{d_j}}.
\end{eqnarray*}
For $\|\zeta\|\geq R$, with $R>0$ large enough, 
one has, since $\deg q_j<d_j $, $j=1,\dots,n$, that
\begin{eqnarray*}
\sum\limits_{j=1}^n{\vert q_j(\zeta)\vert \over \|\zeta\|^{d_j}}<{\tilde c\over 2}.
\end{eqnarray*}
Therefore, for $\|\zeta\|\geq R$, we have 
\begin{eqnarray*}
\sum\limits_{j=1}^n{\vert P_j(\zeta)\vert \over \|\zeta\|^{d_j}}
\geq {\tilde c \over 2}.
\end{eqnarray*}
The last inequality implies {\sl ii}) with 
some constant $c=c(R)$. $\diamondsuit $\\
\vskip 2mm
Note that, if $P$ is a polynomial map from ${\bf C}^n$ to 
${\bf C}^n$, the fact that 
$$
\lim\limits_{\|\zeta\|\mapsto +\infty }\|P(\zeta)\|=+\infty 
$$
(which means just that the map is a proper polynomial 
map in the topological sense) 
does not imply the strong properness condition $(2.1)$. 
For example, if $n=2$, the polynomial map $(X_1X_2,\; (X_1+1)(X_2+1))$ 
is proper, but does not satisfy (2.1) since there  are two common 
zeroes at infinity. 

In order to weaken condition (2.1), we introduce the following concept: 
\begin{definition}
Let $(P_1, \dots , P_n)$ be a polynomial map from ${\bf C}^n$ to 
${\bf C}^n$ and $(\delta _1, \dots , \delta _n)$ 
be a set of strictly positive rational numbers with 
$0<\delta _j\leq \deg P_j $ for any $j$. Then we say that 
$(P_1, \dots , P_n)$ is $(\delta_1, \dots , \delta_n)$- 
proper if and only if there exist $c>0, R>0$ such that, 
for any $\zeta\in {\bf C}^n$ such that $\|{\rm Re}\,\zeta\|\geq R$,  
\begin{eqnarray}
\sum\limits_{j=1}^n{\vert P_j(\zeta)\vert 
\over (1+\| \zeta\| ^2)^{{\delta _j\over 2}}}\geq c.
\end{eqnarray}
\end{definition}
\begin{example} {\rm When $n=2$, the polynomial map $(X_1X_2,(X_1+1)(X_2+1))$ 
is $(1,1)$-proper.}
\end{example}
\begin {remark} {\rm We may extend this notion to the case when the 
$\delta_j$ are rational numbers with the sole conditions 
$\delta_j\leq \deg P_j$. In this setting, a polynomial map which is 
$(\delta_1,\dots,\delta_n)$-proper is not necessarily proper in the 
topological sense.}
\end{remark} 
Let us now formulate the toric analogue of the 
Proposition 2.1. 
\begin{prop}
Let $F_1,\dots F_n $ be n Laurent polynomials with  
Newton polyhedra $\Delta_1,\dots \Delta_n$. 
The following two assertions are equivalent:\\
i) $F_1,\dots F_n $ satisfy the Bernstein conditions (1.3)\\
ii) There exist strictly positive constants $R, c$ such that, 
for any $\zeta \in {\bf C}^n$, with $\|{\rm Re}\, \zeta\|\geq R$,   
\begin{eqnarray}
\sum\limits_{j=1}^n {\vert F_j(e^{\zeta_1}, \dots, e^{\zeta_n}) 
\vert \over e^{H_{\Delta _j}({\rm Re}\, \zeta)}}\geq c\,. 
\end{eqnarray}
\end{prop}
{\bf Proof.}
We first prove that $(i)$ implies $( ii)$. Let us 
assume that $(F_1,\dots,F_n)$ satisfy the Bernstein conditions (1.3). 
In order to prove $(ii)$,  
it is enough to show that one can find 
a conic open sector $S_u$ in ${\bf R}^n$ 
 containing $-u$ and strictly positive 
constants $R_u,\, c_u$, such that, for any $\zeta\in 
{\bf C}^n$ with ${\rm Re}\,\zeta \in S_u$ and $\|{\rm Re}\, \zeta\| 
\geq R_u$, one has
\begin{eqnarray}  
\sum\limits_{j=1}^m 
{\vert F_j(e^{\zeta_1}, \dots , e^{\zeta_n})
\vert \over e^{H_{\Delta _j} ({\rm Re}\, \zeta)}}\geq c_u\,.
\end{eqnarray}
Then, if one can do so for each $u$, 
the existence of positive constants $R$ and $c$ will follow from 
a compactness argument. \\
Applying in the $\zeta$-space a change of coordinates 
$\zeta^{\prime }=A\zeta$, $A\in GL_n({\bf Z})$, we may assume that 
$u=(1,0,\dots ,0)=e_1$. Let us write, for $j=1,\dots,n$, 
\begin{eqnarray}
 F_j(e^{\zeta_1}, \dots , e^{\zeta_n})=e^{k_j\zeta_1}
f_j(e^{\zeta_2}, \dots , e^{\zeta_n})+\tilde F_j(e^{\zeta_1}, \dots , 
e^{\zeta_n})\,,\ 
\end{eqnarray}
where the support of $\tilde F_j$ is included in $\{x_1>k_j\}$.
As noticed by Kazarnovskii \cite{ka1:gnus} (see also \cite{pel:gnus}, 
section 2, from which we inspired ourselves here),   
the fact that Bernstein conditions (1.3) are satisfied for 
 $(F_1, \dots , F_n)$ is equivalent to the following fact: 
for any set of respective faces 
$(\gamma_1,\dots ,\gamma_n)$ of the Newton polyhedra 
$\Delta_1,\dots,\Delta_n$ of $F_1,\dots,F_n$, 
there exists $\epsilon(\gamma_1,\dots,\gamma_n) >0$ 
such that, for any $(\zeta_1,\dots,\zeta_n)\in {\bf C}^n$,  
\begin{eqnarray*}
\sum\limits_{j=1}^n{\vert F_j^{\gamma_j}(e^{\zeta_1}, \dots , 
e^{\zeta_n})
\vert \over e^{H_{\gamma_j}({\rm Re}\, \zeta_1,
\dots, {\rm Re}\,\zeta_n)}}
\geq \epsilon(\gamma_1,\dots,\gamma_n)\,,
\end{eqnarray*}
where, for each $j=1,\dots,n$, 
$F_j^{\gamma_j}$ denotes 
the part obtained from $F_j$ by keeping only monomials corresponding 
to points in $\gamma_j$ and deleting all the others. It is clear 
that whenever $\delta_j$ denotes the Newton polyhedron of $f_j$ 
(considered as a Laurent polynomial in $n-1$ variables  
with support in the subspace $e_1^{\perp}$ of ${\bf R}^n$), the 
convex sets $\tilde\delta_j=
\delta_j+k_j e_1$, $j=1,\dots,n$, are respective 
faces of $\Delta_1,\dots,\Delta_n$. Therefore, one has, for some 
$\epsilon>0$, for $(\zeta_1,\dots,\zeta_n)\in {\bf C}^n$,  
\begin{eqnarray}
\sum\limits_{j=1}^n{\vert e^{k_j\zeta_1}f_j(e^{\zeta_2}, \dots , 
e^{\zeta_n})
\vert \over e^{H_{\tilde \delta_j}({\rm Re}\, \zeta_1,
\dots, {\rm Re}\,\zeta_n)}}
\geq \epsilon\,.
\end{eqnarray}
Since the support of $\tilde F_j$ in (2.5) is included in $\{x_1>k_j\}$, there
exists $\rho>0$, such that, for any $\zeta
=(\zeta_1,\dots,\zeta_n)$ such that ${\rm Re}\, (\zeta_1)<0$ and 
$|{\rm Re}\, \zeta_j|\leq \rho |{\rm Re}\, \zeta_1|$ for 
$j=2,\dots,n$, one has 
\begin{eqnarray}
H_{\tilde \delta_j}({\rm Re}\, \zeta)=H_{\Delta_j}({\rm Re}\, \zeta),\ 
j=1,\dots,n\,.
\end{eqnarray}
On the other hand, if $\rho$ is small enough, then there exists 
$R>0$ such that for any $\zeta\in {\bf C}^n$ such that ${\rm Re}\, \zeta_1
\leq -R$ and $|{\rm Re}\, \zeta_j|\leq \rho|{\rm Re}\, \zeta_1|$ for 
$j=2,\dots,n$, one has
\begin{eqnarray}
\sum\limits_{j=1}^n{\vert \tilde F_j(e^{\zeta_1}, \dots , 
e^{\zeta_n})
\vert \over e^{H_{\tilde \delta_j}({\rm Re}\, \zeta_1,
\dots, {\rm Re}\,\zeta_n)}}
<{\epsilon\over 2}\,.
\end{eqnarray}
From (2.6), (2.7) and (2.8), we get that for $\zeta$ in the conic sector 
$$
S_u:=
\{{\rm Re}\, \zeta_1<0,\ |{\rm Re}\, \zeta_j|<\rho |{\rm Re}\, \zeta_1|,\ 
j=2,\dots,n\}\,,  
$$
the inequality (2.4) is valid for $\|{\rm Re}\, \zeta\|\geq R=R_u$ and 
$c_u=\epsilon/2$. This shows that $(ii)$ holds for the system 
$(F_1,\dots,F_n)$. \\
In order to prove the converse direction we will  construct a globallly defined real analytic 
function that is not vanishing in ${\cal X}\setminus {\bf T}$. This is done as follows:\\
For each $j\in \{1,\dots ,n\}$ choose $n$ 
Laurent polynomials $(G_1^{(j)}, \dots , G_n^{(j)})$ with 
Newton polyhedron $\Delta_j$ such that the system 
$(G_1^{(j)}, \dots , G_n^{(j)})$ satisfies the Bernstein conditions 
(1.3). It follows from the fact that  $(i)$ implies $(ii)$ that, for some convenient constants $C_j\geq c_j>0,\, 
R_j>0$, one 
has, for any $\zeta \in {\bf C}^n$ with $\|{\rm Re}\, \zeta\|\geq R_j$,  
\begin{eqnarray*}
c_j e^{H_{\Delta _j}({\rm Re}\, \zeta)} 
\leq \sum\limits_{k=1}^n\vert G_k^{(j)}
(e^{\zeta_1},\dots , e^{\zeta_n})\vert 
\leq C_j e^{H_{\Delta _j}({\rm Re}\, \zeta)}\,. 
\end{eqnarray*}
Consider now on the torus ${\bf T}^n $ the real analytic function 
$$ 
\zeta\mapsto \phi (\zeta):=
\sum\limits_{j=1}^n{\vert F_j(\zeta)\vert ^2
\over \sum\limits_{k=1}^n\vert G_k^{(j)}(\zeta)\vert ^2}. 
$$
Let ${\cal X}$ be any toric variety associated to a simple 
refinement of the fan which corresponds to 
$\Delta_1+\dots +\Delta_n $. The Laurent polynomials 
$(G_1^{(j)}, \dots , G_n^{(j)})$ induce effective Cartier divisors 
$({\cal D}_1^{(j)}, \dots , {\cal D}_n^{(j)})$ on ${\cal X}$, 
namely 
$$
{\cal D}_k^{(j)}={\rm div}(G_k^{(j)})+E(\Delta_j)\,, \ 1\leq j,k\leq n,
$$
where $E(\Delta_j)$ is the ${\bf T}$-Cartier divisor on ${\cal X}$ 
corresponding to $\Delta_j$ (it is well defined, since 
${\cal X}$ corresponds to a fan which 
is compatible with $\Delta _j$). The fact that 
the system $(G_1^{(j)},\dots,G_n^{(j)})$ obeys Bernstein conditions is 
equivalent (see for example \cite {fu:gnus}) to    
\begin{eqnarray*}
\vert {\cal D}_1^{(j)}\vert \cap \dots \cap \vert {\cal D}_n^{(j)}\vert 
={\cal L}_j\subset {\bf T}^n.
\end{eqnarray*}
For homogeneity reasons, the function
$$
\zeta \mapsto \phi (\zeta_1, \dots , \zeta_n) 
$$
extends from $ {\bf T}^n\setminus \bigcup\limits_{j=1}^n{\cal L}_j$ 
to a function $\tilde \phi $ defined globally as a 
real analytic function on ${\cal X}\setminus 
 \bigcup\limits_{j=1}^n{\cal L}_j$. \\
  Now we are ready to complete the proof of the final step. Assume that  $(F_1,\dots,F_n)$ 
satisfies $( ii)$. For 
$\vert \zeta_1\vert +\dots +\vert \zeta_n\vert 
+{1\over \vert \zeta_1\vert}+\dots +{1\over \vert \zeta_n\vert}$ 
large enough, we have, for some constants $0<\tilde c<\tilde C<\infty$,  
\begin{eqnarray*}
\tilde c \leq \vert \tilde \phi (\zeta_1, \dots , \zeta_n)\vert =
\vert \phi (\zeta_1, \dots , \zeta_n)\vert\leq \tilde C.
\end{eqnarray*}
Therefore $\tilde \phi $ does not vanish on 
${\cal X}\setminus {\bf T}^n$, 
which implies that the effective Cartier divisors 
${\cal D}_j$ induced by the $F_j$ on ${\cal X}$ by 
$$
{\cal D}_j={\rm div}\, (F_j)+E(\Delta_j)
$$
are such that
\begin{eqnarray*}
\vert {\cal D}_1\vert \cap \dots \cap \vert {\cal D}_n\vert \subset {\bf T}^n.
\end{eqnarray*}
This is equivalent to say that the Bernstein conditions are 
fullfilled for the system $(F_1, \dots , F_n)$. 
$\quad\diamondsuit$

In order to weaken the 
properness condition (2.2), we introduce 
the toric analogue of Definition 2.1. 
\begin{definition}
Let $(F_1, \dots , F_n)$ be a system of Laurent polynomials 
in $n$ variables, with Newton polyhedra 
$\Delta_1,\dots,\Delta_n$, and $(\delta _1, \dots , \delta _n)$ 
be a collection of closed convex polyhedra 
with vertices in ${\bf Q}^n$, with $\delta _j \subset 
\Delta_j$, $j=1,\dots,n$. 
Then, we say that $(F_1,\dots ,F_n)$ is 
$(\delta_1, \dots , \delta_n)$-proper 
if and only if there exist $c>0, R>0$ such that, for any $\zeta\in 
{\bf C}^n$ such that $\|{\rm Re}\, \zeta\|\geq R$, 
\begin{eqnarray}
\sum\limits_{j=1}^n
{\vert F_j(e^{\zeta_1}, \dots ,e^{\zeta_n})\vert 
\over e^{H_{\delta _j}({\rm Re}\,\zeta)}}\geq c.
\end{eqnarray}
\end{definition}
\begin{example}
{\rm Let $n=2$ and 
\begin{eqnarray*}
F_1&=&X_1^2X_2^2+X_1^2X_2^{-2}+\alpha_1X_1X_2+\beta_1X_1^{-2}
X_2{-2}+\gamma_1X_2^2X_1^{-2}+
\delta_1X_1^{-2}X_2^{-2}\\
F_2&=&X_1^2X_2^2+X_1^2X_2^{-2}+\alpha_2X_1X_2+\beta_2X_1^{-2}
X_2{-2}+\gamma_1X_2^2X_1^{-2}+\delta_2X_1^{-2}X_2^{-2},
\end{eqnarray*}
with the conditions
\begin{eqnarray*}
\gamma _1(\delta _1-\delta _2)-\delta _1(\gamma _1-\gamma _2)&\not =&0\\
(\alpha _1-\alpha _2)^2-(\beta _1-\beta _2)^2&\not=& 0.
\end{eqnarray*}
Then $(F_1,F_2)$ is $(\delta, \delta)$-proper, where 
$$
\delta={\overline {{\rm conv}\, \{(-2, -2),(2,2),(1,1),(1,-1)\}}}.
$$
In fact, 
it is enough to notice that $(F_1-F_2, F_1)$ 
satisfy the Bernstein conditions and have as respective Newton polyhedra
$\delta$ and $[-2,2]\times [-2,2]$, so that by Proposition 2.2, one has, 
for $\|({\rm Re}\, \zeta_1,{\rm Re}\, \zeta_2)\|\geq R>0$, 
\begin{eqnarray*}
{\vert (F_1-F_2)(e^{\zeta_1}, e^{\zeta_2})\vert 
\over e^{H_{\delta }({\rm Re}\, \zeta_1,{\rm Re}\, \zeta_2)}}+
{\vert F_1(e^{\zeta_1}, e^{\zeta_2})
\vert \over e^{H_{[-2,2]^2}({\rm Re}\, \zeta_1,{\rm Re}\, \zeta_2)}}\geq c\,,
\end{eqnarray*}
which implies, for such $\zeta$,}
\begin{eqnarray*}
{\vert F_1(e^{\zeta_1}, e^{\zeta_2})\vert 
\over e^{H_{\delta }({\rm Re}\, \zeta_1,{\rm Re}\, \zeta_2)}}+
{\vert F_2(e^{\zeta_1}, e^{\zeta_2})
\vert \over e^{H_{\delta}({\rm Re}\, \zeta_1,{\rm Re}\, \zeta_2)
}}\geq {c\over 2}.
\end{eqnarray*}
\end{example}
\section{Proof of the Vanishing Theorems}
\setcounter{equation}{0}
\subsection{ The case of the projective space ${\bf P}^n$}
Our basic tool will be  
multidimensional residue theory through an approach based on the use of  
Bochner-Martinelli integral representation formulaes. 
Let us recall here some well known facts. Let $P_1,\dots,P_n$ 
be $n$ polynomials in $n$ variables defining a discrete (hence 
finite) variety in ${\bf C}^n$. It is shown in \cite{pty:gnus} 
that if $\alpha \in \{P_1=\dots =P_n=0\}$ and 
$\varphi\in {\cal D}({\bf C}^n)$ is such that $\varphi\equiv 1$ 
in a neighborhood of $\alpha$ and $\varphi\equiv 0$ in a neighborhood 
of any point in $\{P_1=\dots=P_n=0\}\setminus \{\alpha\}$,  
we have 
\begin{eqnarray}
&{\rm Res}_{\alpha }[{Qd\zeta \over P_1\dots P_n}]=\nonumber \\
&={(-1)^{{n(n-1)\over 2}}(n-1)!
\over ( 2\pi i)^n}\lim\limits_{\epsilon \mapsto 0}\ 
{1\over \epsilon ^n}
\int\limits_{\| P\|^2=\epsilon }
Q\Big(\sum\limits_{k=1}^n
(-1)^{k-1}\overline {P_k}\bigwedge\limits_{{l=1}\atop{l\not=k}}^n 
\overline {\partial P_l} \Big)\wedge \varphi d\zeta\\
&={(-1)^{{n(n-1)\over 2}}(n-1)! \over ( 2\pi i)^n}
\lim\limits_{\epsilon \mapsto 0}
\int\limits_{\| P\|^2=\epsilon }
{Q\Big (\sum\limits_{k=1}^n
(-1)^{k-1}\overline {P_k}\bigwedge\limits_{{l=1}\atop{l\not=k}}^n 
\overline {\partial P_l}\Big ) \wedge \varphi d\zeta \over \| P\|^{2n}},
\end{eqnarray}
where as usual $\|P\|^2=\vert P_1\vert ^2+\dots +\vert P_n\vert^2 $. 
Using Stokes 's theorem and observing that the form 
$$ 
{Q\Big(\sum\limits_{k=1}^n
(-1)^{k-1}\overline {P_k}\bigwedge\limits_{{l=1}\atop{l\not=k}}^n 
\overline {\partial P_l}\Big) \wedge \varphi d\zeta \over \| P\|^{2n}} 
$$ 
is closed in a punctured neighborhood  
$U_\alpha \setminus \{\alpha\} $, we get from (3.2) 
that if $U_\alpha$ is small enough and with piecewise 
smooth boundary, 
\begin{eqnarray}
{\rm Res}_{\alpha }[{Qd\zeta \over P_1\dots P_n}]
={(-1)^{{n(n-1)\over 2}}(n-1)! \over ( 2\pi i)^n}
\int\limits_{\partial U_{\alpha}}
{Q\Big (\sum\limits_{k=1}^n
(-1)^{k-1}\overline {P_k}\bigwedge\limits_{{l=1}\atop{l\not=k}}^n 
\overline {\partial P_l}\Big ) \wedge d\zeta \over \| P\|^{2n}}\,. 
\nonumber
\end{eqnarray}
Therefore, if $U$ is any bounded open set with smooth boundary containing 
in its interior the set $V(P):=\{P_1=\dots=P_n=0\}$, we have 
\begin{eqnarray} 
&{\rm Res}\left [
 \begin{array}{ccccc}
 Q(X_1,\dots, X_n)dX\\
 P_1,\dots ,P_n\\
 \end{array}
 \right ]=\nonumber \\
&={(-1)^{{n(n-1)\over 2}}(n-1)! \over ( 2\pi i)^n}
\int\limits_{\partial U}
{Q\Big (\sum\limits_{k=1}^n
(-1)^{k-1}\overline {P_k}\bigwedge\limits_{{l=1}\atop{l\not=k}}^n 
\overline {\partial P_l}\Big ) \wedge d\zeta \over \| P\|^{2n}}\,.
\end{eqnarray} 
We can rewrite $(3.3) $ as follows: if 
\begin{eqnarray*}
s_0=({\overline {P_1} 
\over \|P\|^2},\dots ,{\overline {P_n} \over\| P\|^2})=(s_{01},\dots,s_{0n}),
\end{eqnarray*}
then 
\begin{eqnarray*}
{\rm Res}\left [
 \begin{array}{ccccc}
 Q(X_1,\dots, X_n)dX\\
 P_1,\dots ,P_n\\
 \end{array}
\right ] 
=\gamma_n 
\int\limits_{\partial U}
Q(\zeta )\left (\sum\limits_{k=1}^n
(-1)^{k-1}s_{0k}ds_{0,[k]}\right )\wedge d\zeta ,
\end{eqnarray*}
where $ds_{0,[k]}:=
\bigwedge \limits_{j\not=k} ds_{0j}$, $k=1,\dots,n$, and 
$$
\gamma_n:= {(-1)^{{n(n-1)\over 2}}(n-1)! \over ( 2\pi i)^n} \,.
$$ 
An homotopy argument shows that one can replace the 
vector-function $s_0$ above by any vector-function $s$, 
which is $C^1 $ in a neighborhood of the $\partial U$ and satisfies
\begin{eqnarray*}
<s(\zeta ),P(\zeta )>=\sum\limits _{k=1}^n 
s_k(\zeta )P_k(\zeta )\equiv 1 ,\; \zeta \in \partial U.
\end{eqnarray*}
Then the global residue 
is given by the generalized Bochner-Martinelli formula
\begin{eqnarray}
{\rm Res}\left [
\begin{array}{ccccc}
Q(X_1\dots ,X_n)dX\\
 P_1,\dots ,P_n\\
\end{array}
 \right ]=\gamma_n  
\int\limits_{\partial U}
Q \left (\sum\limits_{k=1}^n
(-1)^{k-1}s_{k}ds_{[k]}\right )\wedge d\zeta\,.
\end{eqnarray}
At this stage we are ready for the 
\vskip 2mm
\noindent
{\bf Proof of the Theorem 1.1.}  
\vskip 2mm
\noindent
The first point is to notice that one can assume that the $\delta_j$, 
$j=1,\dots,n$, are strictly positive integers. In order to do so, it 
is enough to use the compatibility of the residue calculus with the 
change of basis (see for example \cite {lip:gnus}, section 2, prop. 2.3), 
which asserts that, for any $N\in {\bf N}^*$, 
\begin{eqnarray}
&{\rm Res}\, \left[ \matrix {Q(X) dX \cr 
P_1(X),\dots,P_n(X)}\right]=\nonumber \\
&={\rm Res}\, \left[ \matrix {Q(X_1^N,\dots, X_n^N) (X_1\cdots X_n)^{N-1} 
dX \cr
P_1(X_1^N,\dots,X_n^N),\dots,P_n(X_1^N,\dots,X_n^N)}\right]\,,
\end{eqnarray}
Let $N$ be a common denominator for the rational numbers $\delta_j$, 
$j=1,\dots,n$; then the polynomials 
$$
\widetilde P_j(X)=P_j(X_1^N,\dots,X_n^N)\,,\ j=1,\dots,n,
$$
have respective degrees $N \deg P_j$, $j=1,\dots,n$, and satisfy (2.2) 
with $\tilde\delta_j=N\delta_j \in {\bf N}^*$. If we assume that our result 
holds when the $\delta_j$ are integers, we get that the residue symbol 
(3.5) is zero when 
$$
N \deg Q +n(N-1) \leq N(\delta_1+\cdots+\delta_n)-n-1\,,
$$  
that is 
$$
\deg Q\leq \delta_1+\cdots+\delta_n-n-{1\over N}\,.
$$
Therefore, we have (1.8) whenever 
$$
\deg Q \leq \delta_1+\cdots+\delta_n-n-1
$$
as we want. We will assume from now on that $\delta_j\in {\bf N}^*$ 
for any $j\in \{1,\dots,n\}$. 
\vskip 2mm
\noindent
Let us denote $D_j:=\deg P_j$ and pick 
an integer $M$ large enough, so that 
$$ 
\delta _k+M-D_k >0, \; \forall k \in \{1, \dots, n\}.
$$
Let $R$ as in the hypothesis of Theorem 1.1, so that, in particular, the 
open ball $B(0,R)$ contains $V(P)$. Let 
us  define 
the vector function $s=s^{\delta,M}$ in ${\bf C}^n \setminus V(P)$ as follows
\begin{eqnarray*}
s^{\delta,M}
(\zeta)=
{1\over \|P(\zeta)\|^2_{\delta,M}}
\Big({\overline {P_1(\zeta)} \over (1+\|\zeta\|^2)^{\delta _1+M}},\dots ,
{\overline {P_n(\zeta)} \over (1+\|\zeta\|^2 )^{\delta_n +M}}\Big)\,,
\end{eqnarray*}
where
$$
\|P(\zeta)\|^2_{\delta,M}:=\sum\limits_{j=1}^n {|P_j(\zeta)|^2 \over 
(1+\|\zeta\|^2)^{M+\delta_j}}\,.
$$
Let 
$$
s_0^{\delta,M}(\zeta):=
\Big({\overline {P_1(\zeta)} \over (1+\|\zeta\|^2)^{\delta _1+M}},\dots ,
{\overline {P_n (\zeta)} \over (1+\|\zeta\|^2 )^{\delta_n +M}}\Big)\,.
$$
Formula (3.4) implies that
\begin{eqnarray}
& {\rm Res}\left [
 \begin{array}{ccccc}
 Q(X)dX\\
 P_1,\dots ,P_n \\
 \end{array}
 \right ]=\nonumber \\ 
&=\gamma _n 
\int\limits_{\|\zeta\|=R}
Q\Big (\sum\limits_{k=1}^n
(-1)^{k-1} s_{k}^{\delta,M }ds_{[k]}^{\delta,M}\Big)\wedge d\zeta=
\nonumber\\
&=\gamma _n  
\int\limits_{\|\zeta\|=R}
\|P\|^{-2n}_{\delta,M} Q\Big (\sum\limits_{k=1}^n
(-1)^{k-1} s_{0k}^{\delta,M }
ds_{0,[k]}^{\delta,M}\Big)\wedge d\zeta=\nonumber 
\\
&=\gamma _n  
\left[\int\limits_{\|\zeta\|=R}
\|P\|^{2(\lambda-n)}_{\delta,M} Q\Big (\sum\limits_{k=1}^n
(-1)^{k-1} s_{0k}^{\delta,M}ds_{0,[k]}^{\delta,M}\Big)\wedge d\zeta
\right]_{\lambda=0}\,.
\end{eqnarray}
For $\lambda$ fixed  with ${\rm Re}\,\lambda >>1 $, let us 
express in homogeneous coordinates 
$\tilde \zeta :=(\zeta_0,\dots , \zeta_n)$ the differential form 
$$ 
\|P\|^{2(\lambda-n)}_{\delta,M} Q\Big (\sum\limits_{k=1}^n
(-1)^{k-1} s_{0k}^{\delta,M }ds_{0,[k]}^{\delta,M}\Big)\wedge d\zeta\,.
$$ 
This leads to a differential $(n,n-1)$ form in ${\bf P}^n$ (depending on the 
complex parameter $\lambda$), which will be denoted as 
$\Theta_{P, Q,\lambda}^{\delta,M}$. Since 
$$
\overline\partial
\left[\|P\|^{2(\lambda-n)}_{\delta,M} Q\Big (\sum\limits_{k=1}^n
(-1)^{k-1} s_{0k}^{\delta,M}ds_{0,[k]}^{\delta,M}\Big)\wedge d\zeta\right]
=n\lambda 
\|P\|^{2(\lambda-n)}_{\delta,M}Q 
\Big(\bigwedge\limits_{k=1}^n \overline\partial 
s_{0k}^{\delta,M}\Big)\wedge d\zeta\,,
$$
we have, if the action of the 
$\overline\partial$ operator is now considered on the projective
differential forms (expressed in homogeneous coordinates), 
\begin{eqnarray}
&\overline\partial \Theta_{P, Q,\lambda}^{\delta,M}=\nonumber \\
&=n\lambda \Big({\|{\cal P}(\tilde \zeta)\|_{\delta,M}\over \|\tilde 
\zeta\|^{|\delta|+M}}\Big)^{2(\lambda-n)}
\zeta_0^{-\deg Q-n-1} {\cal Q}(\tilde \zeta)
A_{P,Q}^{\delta,M}
({\zeta_1\over \zeta_0},\dots,{\zeta_n\over \zeta_0})\wedge 
\Omega(\tilde \zeta)\,,
\end{eqnarray}
where $|\delta|=\delta_1+\cdots+\delta_n$, $\delta_{[j]}=|\delta|-\delta_j$ 
for $j=1,\dots,n$, $\Omega$ is the Euler form, 
$$
A_{P,Q}^{\delta,M}:=
\bigwedge\limits_{k=1}^n \overline\partial 
s_{0k}^{\delta,M}\,,
$$
and
$$
\|{\cal P}(\tilde \zeta)\|^2_{\delta,M}:=\sum\limits_{j=1}^n 
\vert {\cal P}_j\vert ^2\vert \zeta_0\vert ^{2(\delta_j-D_j+M)}
 \|\tilde \zeta \|^{2\delta _{[j]}}\,,
$$ 
${\cal P}_1,\dots, {\cal P}_n,{\cal Q}$, being the respective 
homogeneizations of $P_1,\dots,P_n,Q$; the norm $\|\tilde \zeta\|$ 
is the Euclidean norm in ${\bf C}^{n+1}$.  
Since 
$$
s_{0k}^{\delta,M}({\zeta_1\over \zeta_0},\dots,{\zeta_n\over \zeta_0})
=|\zeta_0|^{2(\delta_k+M)} {\overline {{\cal P}_k(\tilde \zeta)} 
\overline {\zeta_0}^{-D_k} \over \|\tilde \zeta\|^{2(\delta_k+M)}}\,,
$$
one has
$$
A_{P,Q}^{\delta,M}({\zeta_1\over \zeta_0},\dots, 
{\zeta_n\over \zeta_0})=
\zeta_0^{nM+|\delta|} 
\bigwedge\limits_{k=1}^n
\overline\partial \left[{\overline \zeta_0^{\delta_k+M-D_k}
\overline {{\cal P}_k} \over \|\tilde \zeta\|^{2(M+\delta_k)}}\right]\,.
$$   
In the same vein, we have 
\begin{eqnarray*}
&\Theta_{P, Q,\lambda}^{\delta,M}(\tilde \zeta)=
\big({\|{\cal P}\|_{\delta,M}\over 
\|\tilde \zeta\|^{|\delta|+M}}\big)^{2(\lambda-n)} 
\zeta_0^{nM+|\delta|-n-1}
{\cal Q}(\tilde \zeta) \times \\
&\times \Bigg(\sum\limits_{k=1}^n 
(-1)^{k-1} {\overline \zeta_0^{\delta_k+M-D_k}
\overline {{\cal P}_k} \over \|\tilde \zeta\|^{2 (\delta_{[k]}+M)}}
\bigwedge\limits_{{l=1}\atop {l\not=k}}^k 
\overline\partial \left[{\overline \zeta_0^{\delta_l+M-D_l}
\overline {{\cal P}_l}\over \|\tilde \zeta\|^{2 (\delta_{[l]}+M)}}
\right]\Bigg)
\wedge \Omega\,.
\end{eqnarray*}
This shows (as a consequence of Atiyah's theorem \cite{at:gnus}) that the map 
$$
\lambda \mapsto \Theta_{P, Q,\lambda}^{\delta,M}
$$
can be considered as a meromorphic map with values in the 
space of $(n,n-1)$ currents in ${\bf P}^n({\bf C})$.  
\vskip 1mm
\noindent
We now consider the complement in ${\bf P}^n({\bf C})$ of $B(0,R)$ 
as a $2n$-chain $\Sigma$ in ${\bf P}^n$ (with smooth boundary). 
One has, for ${\rm Re}\, \lambda>>1$, using Stokes's theorem 
$$
\int_{\partial \Sigma} \Theta_{P, Q,\lambda}^{\delta,M}
=\int_{\Sigma} \overline \partial [\Theta_{P, Q,\lambda}^{\delta,M}]\,.
$$  
Therefore, one can rewrite (3.6) as 
\begin{eqnarray}
{\rm Res}\left[
\begin{array}{ccccc}
Q(X)dX\\
P_1,\dots ,P_n \\
\end{array}
\right]=-\gamma_n 
\left[\int_{\partial \Sigma} \Theta_{P, Q,\lambda}^{\delta,M}
\right]_{\lambda=0}
=-\gamma_n\left[\int_{\Sigma} \overline\partial 
\Big(\Theta_{P, Q,\lambda}^{\delta,M}\Big)
\right]_{\lambda=0}\,
\end{eqnarray} 
(the total sum of residues in ${\bf C}^n$ equals the opposite of 
the ``residue'' at infinity). In order to compute this residue at infinity 
(and to prove that it vanishes in the situation we are dealing with), 
we localize the problem and look 
at the analytic continuation up to the origin of the meromorphic function 
\begin{eqnarray}
\lambda \mapsto \int_{\Sigma} \varphi \overline\partial 
\Big(\Theta_{P, Q,\lambda}^{\delta,M}\Big)\,,
\end{eqnarray}
when $\varphi$ is an element in ${\cal D}({\bf P}^n({\bf C}))$ with 
support contained in a neighborhood $V$ of some point $x$ at infinity in 
${\bf P}^n({\bf C})$ (these are the only 
interesting points, since if the support of $\varphi$ does not intersect the 
hyperplane at infinity, then (3.9) is an entire function which vanishes 
at $\lambda=0$). We may suppose that the local coordinates in $V$ 
are 
$\xi:=({\zeta_0\over \zeta_1},\dots,{\zeta_n\over \zeta_1})$ (for example). 
Let 
$$
f_j(\zeta)={{\cal P}_j(\tilde \zeta) \zeta_0^{\delta_j+M-D_j}\over 
\zeta_1^{M+\delta_j}},\ j=1,\dots,n,
$$
expressed in the local coordinates $\xi$. Let us introduce a resolution of 
singularities $({\cal X},\pi)$ for the hypersurface $\{f_1\cdots f_n=0\}$ 
over $V$ (schrinking $V$ about the point $x$ if necessary). 
Then, in a local 
chart $\omega$ on ${\cal X}$ with coordinates $w$ centered at the origin, all 
functions $\pi ^{*}(f_j)$ are, up to invertible holomorphic functions,  
monomials in $w $; that is
$$
\pi ^{*}(f_j)(w)=u_j(w) w_1^{\theta_{j1}}\cdots w_n^{\theta_{jn}},\ 
\theta_{jk}\in {\bf N},\ u_j \ {\rm invertible\ in}\ \omega\,.
$$
Note that 
$$
\pi^* \Big[{\zeta_0\over \zeta_1}\Big](w)=u_0(w) w_1^{\theta_{01}}\cdots 
w_n^{\theta_{0n}},\ \theta_{0k}\in {\bf N},\ u_0\ 
{\rm invertible\ in}\ \omega\,,
$$
since $\delta_j+M-D_j>0$ for at least one $j$ (in fact for any $j$). 
However this is not enough. 
Using the ideas of A. Varchenko \cite{va:gnus} 
and A. Khovanskii \cite{kh1:gnus}, we introduce, 
above each such local chart $\omega$, a toroidal 
manifold $\widetilde{\cal X}$ and a proper holomorphic map 
$\tilde \pi: \widetilde {\cal X}\mapsto \omega$ (wich is locally 
a biholomorphism between $\widetilde {\cal X}\setminus \tilde\pi^* 
\{w_1\cdots w_n=0\}$ and $\omega \setminus \{w_1\cdots w_n=0\}$), such that, 
on each local chart $\tilde \omega$ on $\widetilde {\cal X}$ (with 
local coordinates $(t_1,\dots,t_n)$), one has 
$$
\tilde\pi^*\pi^* (f_j)(t_1,\dots,t_n)=\tilde u_j(t) t_1^{\tilde \theta_{j1}}
\cdots t_n^{\tilde \theta_{jn}}=\tilde u_j(t) \tilde m_j(t)
$$
and one of the $\tilde m_j$, $j=1,\dots,n$, let say $\tilde m$, 
divides $\tilde m_1,\dots \tilde m_n$.  
Namely, the manifold $\widetilde {\cal X}$ is the smooth 
toric variety attached 
to a simple refinement of the fan associated with the Newton polyhedron 
$$
\Gamma^+(\theta_1,\dots,\theta_n):= 
\bigcup\limits_{j=1}^n \Big[\theta_j+{\bf N}^n\Big]\,.
$$ 
It arises from glueing together copies 
$(U_J,\pi_J)$ of ${\bf C}^n$ (in correspondence 
with the $n$-dimensional cones of the fan, 
$\pi_J$ being a monoidal transform attached to the skeleton of the cone), 
according to the glueing of the cones along  
their edges. The $1$-dimensional edges of these 
cones are determined as the normal 
directions to the $(n-1)$-dimensional faces of the 
Newton polyhedron $\Gamma ^+(\theta_1,\dots, \theta_n)$, plus a
minimal system of additional directions rational directions in 
$[0,\infty[^n$ (which are just necessary for the fan to be simple).
\vskip 1mm
\noindent
We now come to the crucial point where we use the hypothesis (1.7), which 
tells us that, for 
$R |\zeta_0|\leq (|\zeta_1|^2+\cdots+|\zeta_n|^2)^{1\over 2}$, one 
has
$$
|\zeta_0|^M \|\tilde \zeta\|^{|\delta|}
\leq c \sum\limits_{j=1}^n 
|{\cal P}_j(\tilde \zeta)||\zeta_0|^{M+\delta_j-D_j}
\|\tilde \zeta\|^{\delta_{[j]}}
\leq c_n \|{\cal P}\|_{\delta,M}\,.
$$
This implies that, if 
$$
\tilde\pi^*\pi^* \Big[{\zeta_0\over \zeta_1}\Big](t)
=\tilde u_0(t)t_1^{\tilde \theta_{01}}\cdots t_n^{\tilde\theta_{0n}},\ 
\tilde u_0\ {\rm invertible\ in}\ \tilde\omega\,,
$$
the distinguished monomial $\tilde m$ divides 
$\tilde m_0^M$, where $\tilde m_0:=
t_1^{\tilde \theta_{01}}\cdots t_n^{\tilde\theta_{0n}}$, in $\tilde \omega$. 
Let $\tilde \varphi$ be a test function on ${\cal X}$ with support in 
the local chart $\omega$. 
As one can see it easily, one can write in $\tilde\omega$,  
$$
\tilde\pi^*\bigg[\tilde 
\varphi \bigg[\pi^*\Bigg(\varphi A_{P,Q}\Big({\zeta_1\over \zeta_0},\dots,
{\zeta_n\over \zeta_0}\Big)\Bigg)\bigg]\bigg](t)=
{\tilde m_0^{nM+|\delta|}\over \tilde 
m^n} \Bigg({\overline {\partial \tilde m}\over 
\overline {\tilde m}}\wedge \sigma_1(t) +\tau_1(t)\Bigg)\,,
$$
where $\sigma_1$ and $\tau_1$ (depending on $\varphi$ and $\tilde \varphi$) 
are smooth differential forms with respective type $(n,n-1)$ and $(n,n)$. It 
follows then from (3.7) that
$$
\tilde\pi^*[\tilde 
\varphi [\pi^*\big(\varphi 
\overline\partial \Theta_{P,Q,\lambda}^{\delta,M}\big)]](t)
=\lambda |\tilde m|^{2\lambda} |\xi|^{2\lambda}
{\tilde m_0^{nM+|\delta|-\deg Q-n-1}\over \tilde m^n}
 \Bigg({\overline {\partial \tilde m}\over 
\overline {\tilde m}}\wedge \sigma_2(t) +\tau_2(t)\Bigg)
$$
where $\sigma_2$ and $\tau_2$ (depending on $\varphi$ and $\tilde \varphi$) 
are smooth differential forms with respective type $(n,n-1)$ and $(n,n)$ and 
$\xi$ is a real analytic strictly positive function in $\tilde\omega$. 
Since $|\delta|-\deg Q-n-1\geq 0$ and $\tilde m^n$ divides $\tilde m_0^{nM}$, 
we get immediately that for any test function $\rho$ 
with support in $\tilde \omega$, 
$$
\left[\int \rho(t)\tilde\pi^*[\tilde 
\varphi [\pi^*\big(\overline\partial \Theta_{P,Q,\lambda}^{\delta,M}\big)](t)
\right]_{\lambda=0}=0\,.
$$
Then, the conclusion (1.8) follows from the formula (3.8) and our 
localization and normalized blowing-up process. $\quad\diamondsuit$ 
\begin {remark} 
{\rm 
The fact that $\delta_j>0$ does not play any role in the proof. Therefore, 
theorem 1.1 remains valid when $(P_1,\dots,P_n)$ is 
$(\delta_1,\dots,\delta_n)$-proper, where the $\delta_j$ are rational 
numbers such that $\delta_j\leq D_j$ for any $j=1,\dots,n$ (see Remark 2.1); 
of course, the conclusion of the theorem 
is interesting only in the case when $\delta_1+\cdots+\delta_n 
\geq n+1$.} 
\end{remark}
\subsection{ The toric case}
We begin with a review of some preliminary material 
taken from \cite{ccd:gnus, cd:gnus, cox1:gnus, cox2:gnus, fu:gnus}. \\
A complete toric variety ${\cal X}$ of dimension $n$ 
is determined by a complete fan ${\cal F}$ in an $n$-dimensional 
real vector space $\Lambda_{{\bf R}}$, where $\Lambda$ is a lattice; for the 
sake of simplicity, we will always assume $\Lambda={\bf Z}^n$ and 
$\Lambda_{\bf R}={\bf R}^n$. Taking a suitable refinement of the fan, 
we may assume that this toric variety ${\cal X}$ is also smooth. 
\vskip 1mm
\noindent
We denote as $\Lambda^*\simeq {\bf Z}^n$ the dual lattice. 
The primitive generators of the one dimensinal cones in ${\cal F}$
are denoted by $\eta_1,\dots ,\eta_s$. Each of these vectors $\eta_i$, 
$i=1,\dots,s$, is in correspondence with a torus-invariant irreducible 
Weil divisor ${\cal X}_i$ on ${\cal X}$. The 
$(n-1)$-Chow group ${\bf A}_{n-1}({\cal X})$ on ${\cal X}$ 
is generated by the 
classes $[{\cal X}_i]$, $i=1,\dots,s$, and induces a grading on the 
polynomial algebra ${\bf S}={\bf C}[x_1,\dots,x_s]$, namely 
$$
\deg (x_1^{\alpha_1}\cdots x_s^{\alpha_s}):= 
\Big[\alpha_1 {\cal X}_1+\cdots+\alpha_s 
{\cal X}_s\Big]\in {\bf A}_{n-1}({\cal X})\,.
$$ 
Note that the sequence 
$$ 
0\rightarrow \Lambda ^* \buildrel {\tau}\over {\rightarrow} 
{\bf Z}^s\rightarrow {\bf A}_{n-1}({\cal X})\rightarrow 0\,,
$$
where $\tau(m)=(<m,\eta_1>,\dots,<m,\eta_s>)\in {\bf Z}^s$ is exact 
since any monomial $x^{<m,\eta >}:=x_1^{<m,\eta _1>}\dots x_s^{<m,\eta _s>}$, 
$m\in \Lambda^*$, has degree zero. If $(e_1^*,\dots,e_n^*)$ is 
the canonical basis of $\Lambda^*$ and ${\cal I}$ is an
ordered  subset of 
$\{1,...,s\}$ with cardinal $|{\cal I}|=n$, let say  
${\cal I}=\{i_1,\dots,i_n\}$, 
$1\leq i_1<\cdots<i_n\leq s$, we denote as 
$$
dx_{\cal I}:=\bigwedge\limits_{l=1}^n dx_{i_l},\ 
\widehat {x_{\cal I}}:=\prod\limits_{{k=1}\atop {k\notin {\cal I}}}^s
x_k,\ \det[\eta_{\cal I}]:= \det[<e_k^*,\eta_{i_l}>]_{1\leq k,l\leq n}\,.
$$
The toric Euler form on ${\cal X}$ is the differential form $\Omega$ 
(expressed in homogeneous coordinates $x_1,\dots,x_s$) 
$$
\Omega(x):=
\pm 1 \sum\limits _{\vert {\cal I}\vert =n} \det[\eta _{\cal I}]
\widehat {x_{\cal I}} dx_{\cal I}\,.
$$     
\vskip 2mm
\noindent 
We now consider a system $(F_1,...,F_n)$ of 
Laurent polynomials with respective polyhedra 
$\Delta_1$,...,$\Delta_n$, and a collection $(\delta_1,...,\delta_n)$ of 
rational polyhedra such that $\delta_j\subset \Delta_j$ for any 
$j\in \{1,\dots,n\}$, $\delta_1+\cdots+
\delta_n$ is $n$-dimensional and the hypothesis (1.9) are fulfilled.  
\vskip 2mm
\noindent
Before proceeding any further,
by using the same change of basis (namely replace $X_j$ by $X_j^N$ for 
a convenient choice of $N\in {\bf N}^*$) as we did in (3.5), we can reduce 
ourselves to the situation where all polyhedra $\delta_1$,...,$\delta_n$ 
have their vertices in the lattice $\Lambda={\bf Z}^n$ (originally these 
vertices were assumed to be in ${\bf Q}^n$, therefore it is enough to 
take for $N$ a common denominator of all coordinates of such points).  
\vskip 2mm
\noindent
We fix a polyhedron $\Delta$ with dimension $n$ and vertices in $\Lambda$, 
which contains the origin as an interior point and is such 
that, for any $j\in \{1,\dots,n\}$, the Minkowski sum $\Delta+\delta_j$ 
contains $\Delta_j$. We let 
$$
\widetilde \Delta:=[\Delta+\delta_1+\cdots+\delta_n]+\Delta_1+\cdots+
\Delta_n\,.
$$
We consider as the fan ${\cal F}$ a simple refinement of the fan 
${\cal F}(\widetilde \Delta)$ which corresponds to this polyhedron $\widetilde 
\Delta$  
(see \cite {fu:gnus}); ${\cal X}$ will be from now on the 
toric variety attached to ${\cal F}$. It is 
compatible with $\Delta$, $\delta_j+\Delta$ and $\Delta_j$ for 
any $j$. For any $j=1,\dots,n$, we take 
$n+1$ Laurent polynomials, with convex polyhedron $\Delta+\delta_j$, 
namely $G_0^{(j)},\dots,G_n^{(j)}$, which do not vanish simultaneously in 
${\bf T}^n$ and are such that the system $(G^{(j)}_1,\dots, 
G^{(j)}_n)$ satisfies the Bernstein conditions 
(1.3) (when considered as a system of Laurent polynomials 
with Newton polyhedron 
$\Delta+\delta_j$). Since the fan ${\cal F}$  
is compatible with $\Delta+\delta_j$, 
these Laurent polynomials induce Cartier divisors 
${\cal D}_0^{(j)},\dots,{\cal D}_n^{(j)}$ on ${\cal X}$ such that 
$$
|{\cal D}_0^{(j)}|\cap \dots \cap |{\cal D}_n^{(j)}|=\emptyset\,.
$$  
In particular, the function 
$$
\|G^{(j)}\|^2:= \sum\limits_{k=0}^n |G_k^{(j)}|^2
$$
does not vanish on the torus ${\bf T}^n$. Let, for 
$\zeta\in {\bf T}^n$, 
$$
\|F(\zeta)\|^2_{\delta,\Delta}:=\sum\limits_{j=1}^n 
{|F_j(\zeta)|^2 \over \|G^{(j)}(\zeta)\|^2}\,,
$$ 
$$
s_0^{\delta,\Delta}(\zeta):=
\Bigg({\overline {F_1(\zeta)}\over \|G^{(1)}(\zeta)\|^2},\dots,
{\overline {F_n(\zeta)}\over \|G^{(n)}(\zeta)\|^2}\Bigg)\,,
$$
and, for $\zeta\in {\bf T}^n\setminus V^*(F)$, 
$$
s^{\delta,\Delta}(\zeta):={s_0^{\delta,\Delta}(\zeta)
\over \|F(\zeta)\|^2_{\delta,\Delta}}\,.
$$
Let $\epsilon=\min \{\|\zeta-\zeta'\|;\ \zeta\not=\zeta',\ \zeta,\zeta'\in 
V^*(F)\}$ and 
$$
U:=\bigcup\limits_{\alpha \in V^*(F)} B\big[\alpha, 
{\min\limits_{\alpha \in V^*(F)}(\epsilon,\; 
d(\alpha, {\bf C}^n\setminus {\bf T}^n))\over 2}\big]\, 
$$
where $d$ is the Euclidean distance in ${\bf C}^n$. We can state the 
following 
\begin{lemma}
Let $F=(F_1,\dots ,F_n)$ be 
a system of Laurent polynomials in n variables, with 
respective Newton polyhedra $\Delta_1,\dots,\Delta_n$ 
and polyhedra $\delta_1$, ..., $\delta_n$, $\Delta$ as above. 
Then, for any  Laurent 
momomial $Q(\zeta)=\zeta_1^{\beta_1}\cdots \zeta_n^{\beta_n}=\zeta^\beta $
 one
has
\begin{eqnarray}
&{\rm Res}\,\left [
\begin{array}{ccccc}
Q(X_1,\cdots, X_n)dX\\
F_1,\cdots ,F_n\\
\end{array}
\right ]_{{\bf T}}=\nonumber\\
&=\gamma _n  
\left[\int\limits_{\partial U}
\|F\|^{2(\lambda-n)}_{\delta,\Delta}\zeta^{\beta}
 \Big (\sum\limits_{k=1}^n
(-1)^{k-1} s_{0k}^{\delta,\Delta}ds_{0,[k]}^{\delta,\Delta}\Big)\wedge 
 {d\zeta_1
\over \zeta_1}\wedge \cdots \wedge {d\zeta_n
\over \zeta_n}\right]_{\lambda=0}\,.
\end{eqnarray}
\end{lemma}

\vfill\eject
\noindent
{\bf Proof.} 
\vskip 1mm
\noindent
One has, as in the projective situation (3.6), 
\begin{eqnarray*}
& {\rm Res}\left [
 \begin{array}{ccccc}
 X_1^{\beta_1}\cdots X_n^{\beta_n}\ dX\\
 F_1,\dots ,F_n \\
 \end{array}
 \right ]_{\bf T}=\nonumber \\ 
&=\gamma _n 
\int\limits_{\partial U}
\zeta^{\beta} \Big (\sum\limits_{k=1}^n
(-1)^{k-1} s_{k}^{\delta,\Delta }ds_{[k]}^{\delta,\Delta}\Big)\wedge {d\zeta_1
\over \zeta_1}\wedge \cdots \wedge {d\zeta_n
\over \zeta_n} =
\nonumber\\
&=\gamma _n  
\int\limits_{\partial U}
\|F\|^{-2n}_{\delta,\Delta}
\zeta^{\beta}  \Big (\sum\limits_{k=1}^n
(-1)^{k-1} s_{0k}^{\delta,\Delta}
ds_{0,[k]}^{\delta,\Delta}\Big)\wedge {d\zeta_1
\over \zeta_1}\wedge \cdots \wedge {d\zeta_n
\over \zeta_n} =\nonumber 
\\
&=\gamma _n  
\left[\int\limits_{\partial U}
\|F\|^{2(\lambda-n)}_{\delta,\Delta}\zeta^{\beta}
 \Big (\sum\limits_{k=1}^n
(-1)^{k-1} s_{0k}^{\delta,\Delta}ds_{0,[k]}^{\delta,\Delta}\Big)\wedge 
 {d\zeta_1
\over \zeta_1}\wedge \cdots \wedge {d\zeta_n
\over \zeta_n}\right]_{\lambda=0}\,.
\end{eqnarray*}
This concludes the proof of the lemma. $\diamondsuit $\\
We are now ready to embark ourselves in the 
\vskip 1mm
\noindent
{\bf Proof of Theorem 1.2.} 
We begin with the toric analogue of the standard homogenization 
in the projective spaces.\\
For  $\lambda$ fixed with ${\rm Re}\,\lambda>>1$, we express in homogeneous 
coordinates $(x_1,\dots,x_s)$ the differential form
\begin{eqnarray}
&\Theta_{F,\beta,\lambda}^{\delta,\Delta}:= 
\|F\|^{2(\lambda-n)}_{\delta,\Delta} \zeta^{\beta}
\Big (\sum\limits_{k=1}^n
(-1)^{k-1} s_{0k}^{\delta,\Delta}ds_{0,[k]}^{\delta,\Delta}\Big)\wedge d\zeta
\,;
\end{eqnarray}
taken from $(3.10)$, and where the coordinates $\zeta_j$, $j=1,\dots,n$, in the torus are expressed in 
homogeneous coordinates as 
\begin{eqnarray}
\zeta_j=\prod\limits _{i=1}^s x_i^{\eta _{ij}}
:=x^{<e_j^*,\eta>}\,,\ j=1,\dots,n,
\end{eqnarray}
where, for any $i=1,\dots,s$,  $\eta_{ij}$, $j=1,\dots,n$, are the 
coordinates of the primitive vector $\eta_i$ in the canonical basis 
$(e_1,\dots,e_n)$ of $\Lambda \simeq{\bf Z}^n$. 
In order to do that, we need to introduce the 
${\cal X}$-homogenizations of $F_1,\dots,F_n$, that is 
$$
{\cal F}_j(x_1,\dots,x_s):= 
\Bigg(\prod\limits_{i=1}^s x_i^{\mu_{ij}}
\Bigg) F_j(x^{<e_1^*,\eta>},\dots,x^{<e_n^*,\eta>})\,,
$$
where 
$$
\mu_{ij}:=-\min\limits_{\xi \in \Delta_j} <\xi,\eta_i>,\ i=1,\dots,s;\ 
j=1,\dots,n,
$$
and the ${\cal X}$-homogenizations of the $G^{(j)}_k$, $j=1,\dots,n$, 
$k=0,\dots,n$, namely 
$$
{\cal G}_k^{(j)}(x_1,\dots,x_s)
:=\Bigg(\prod\limits_{i=1}^s 
x_i^{\nu_{ij}}
\Bigg) G_k^{(j)}(x^{<e_1^*,\eta>},\dots,x^{<e_n^*,\eta>})\,,
$$
where 
$$
\nu_{ij}:=-\min\limits_{\xi \in \delta_j+\Delta} <\xi,\eta_i>,\ 
i=1,\dots,s;\ 
j=1,\dots,n\,.
$$
We will also denote 
$$
\|{\cal G}^{(j)}(x_1,\dots,x_s)\|^2:= 
\sum\limits_{k=0}^n |{\cal G}^{(j)}_k(x)|^2,\ j=1,\dots,n\,.
$$
The function 
$$
\zeta \mapsto \|F(\zeta)\|^2_{\delta,\Delta}
$$
on the torus will be extended as the function on ${\cal X}$ which is 
defined in homogeneous coordinates as 
$$
\|{\cal F}(x)\|^2_{\delta,\Delta}:=\sum\limits_{j=1}^n 
\Bigg|\prod\limits_{i=1}^s 
x_i^{\nu_{ij}-\mu_{ij}}\Bigg|^2
\ {|{\cal F}_j(x)|^2 \over \|{\cal G}^{(j)}(x)\|^2}\,.
$$ 
On the other end, one has, for $k=1,\dots,n$,  
$$
s_{0k}^{\delta,\Delta} (x^{<e_1^*,\eta>},\dots,x^{<e_n^*,\eta>})
=\Bigg|\prod\limits_{i=1}^s 
x_i^{\nu_{ik}}\Bigg|^2 
{\Bigg(\prod\limits_{i=1}^s \overline {x_i^{-\mu_{ik}}}
\Bigg) \overline {{\cal F}_k(x)}
\over \|{\cal G}^{(k)}(x)\|^2}\,,
$$
while the differential form 
$$
\Bigg(\prod\limits_{i=1}^s x_i^{<\beta,\eta_i>}
\Bigg)\ {\Omega (x) \over x_1\dots x_s}
$$
on ${\cal X}$ restricts to the torus as $\zeta^\beta 
{d\zeta_1\over \zeta_1}\wedge 
\dots \wedge {d\zeta_n\over \zeta_n}$ (see the proof of 
Proposition 9.5 in \cite{bc:gnus}).
One has 
\begin{eqnarray*} 
&\Theta_{F,\beta,\lambda}^{\delta,\Delta}(x)
=\|{\cal F}(x)\|^{2(\lambda-n)}_{\delta,\Delta}
\prod\limits_{i=1}^s x_i^{<\beta,\eta_i>-1+\sum\limits_{j=1}^n 
\nu_{ij}} \times \\ 
&\times \Bigg( 
\sum\limits_{k=1}^n (-1)^{k-1} 
{\Big[\big(\prod\limits_{i=1}^s \overline {x_i^{\nu_{ik}-
\mu_{ik}}}\big) \overline {{\cal F}_j(x)}\Big]
\over \|{\cal G}^{(k)}(x)\|^2}
\bigwedge\limits_{{l=1}\atop {l\not=k}}^n 
\overline\partial 
\left[ {\Big(\prod\limits_{i=1}^s \overline {x_i^{\nu_{il}
-\mu_{il}}}\Big) \overline {{\cal F}_l(x)} 
\over \|{\cal G}^{(l)}(x)\|^2 }\right]\Bigg)\wedge {\Omega(x)} 
\end{eqnarray*}
and, as in the projective case, 
\begin{eqnarray} 
&\overline \partial \Theta_{F,\beta,\lambda}^{\delta,\Delta}(x)=
\nonumber \\
&=n\lambda
\|{\cal F}(x)\|^{2(\lambda-n)}_{\delta,\Delta}
\Bigg(\prod\limits_{i=1}^s 
x_i^{<\beta,\eta_i>-1+\sum\limits_{j=1}^n 
\nu_{ij}}\Bigg) \bigwedge\limits_{k=1}^n 
\overline\partial 
\left[ {\Big(\prod\limits_{i=1}^s \overline {x_i^{\nu_{ik}
-\mu_{ik}}}\Big) \overline {{\cal F}_k(x)} 
\over \|{\cal G}^{(k)}(x)\|^2 }\right] \wedge \Omega(x)
\,\nonumber\\
& 
\end{eqnarray}
(by the above equalities, we mean that the 
differential forms on ${\cal X}$ which are such defined restrict 
respectively to the 
torus as the differential forms  
$\Theta_{F,\beta,\lambda}^{\delta,\Delta}(\zeta)$ in (3.11) and 
its $\overline \partial$ in $\zeta$). 
\vskip 1mm
\noindent
We need now to interpret our hypothesis (1.9). Since the system of polynomials 
$G_k^{(j)}$, $k=1,...,n$, satisfies the Bernstein hypothesis (1.3), it 
follows (see the argument used in the proof of Proposition 2.2) that 
there exist strictly positive constants $c_j$, $C_j$, such that
\begin{eqnarray}
\forall \zeta \in {\bf C}^n\,,\ 
c_j e^{H_{\Delta+\delta_j}({\rm Re}\, \zeta)}
\leq \|G^{(j)} (e^{\zeta_1},\dots,e^{\zeta_n})\|
\leq C_j  e^{H_{\Delta+\delta_j}({\rm Re}\, \zeta)}\,.
\end{eqnarray}
One has also, for any $\zeta\in {\bf C}^n$ such that 
$\|{\rm Re}\, \zeta\|\geq R$, 
\begin{eqnarray} 
\sum\limits_{j=1}^n {|F_j(e^{\zeta_1},\dots,e^{\zeta_n})|\over 
e^{H_{\delta_j} ({\rm Re}\, \zeta)}}\geq c>0
\end{eqnarray}
and $H_{\delta_j+\Delta}=H_{\delta_j}+H_\Delta$. We also introduce 
$n+1$ Laurent polynomials $H_0,\dots,H_n$, with Newton polyhedron 
$\Delta$, which do not 
vanish simultaneously in 
${\bf T}^n$ and are 
such that the system $(H_1,\dots,H_n)$ satisfies 
the Bernstein hypothesis (1.3) when considered as a system 
of Laurent polynomials with Newton polyhedron $\Delta$, that is such that 
\begin{eqnarray}
c_0\,  e^{H_\Delta ({\rm Re}\, \zeta)}
\leq \|H(e^{\zeta_1},\dots, e^{\zeta_n})\|\leq C_0 \, 
e^{H_\Delta ({\rm Re}\, \zeta)},\ \zeta\in {\bf C}^n 
\end{eqnarray} 
for some strictly positive constants $c_0,C_0$ (where $\|H\|^2:=
|H_0|^2+\cdots+|H_n|^2$). 
It follows from (3.14), (3.15) and (3.16) 
that for any $\zeta\in {\bf T}^n$ such that 
$$
|\zeta_1|+\cdots+|\zeta_n|+{1\over |\zeta_1|}+\cdots 
+{1\over |\zeta_n|}
$$
is large enough, one has 
\begin{eqnarray}
\sum\limits_{j=1}^n {|F_j(\zeta)|\over \|G^{(j)}(\zeta)\|}
\geq {\tilde c\over \|H(\zeta)\|}
\end{eqnarray}
for some $\tilde c>0$. 
If we express the $\zeta_j$ in terms of homogeneous coordinates on the 
toric variety ${\cal X}$ as in (3.12), we may rewrite (3.17) as 
\begin{eqnarray}
\Big|\prod\limits_{i=1}^s x_i^{-\min\limits_{\xi\in {\Delta}}
<\xi,\eta_i>}\Big|\leq {1\over \tilde c} \ {1\over |{\cal H}(x)|}
\Bigg(\sum\limits_{j=1}^n 
\Big|\prod\limits_{i=1}^s x_i^{\nu_{ij}-\mu_{ij}}\Big|\ 
{|{\cal F}_j(x)|\over \|{\cal G}^{(j)}(x)\|}\Bigg)\,,
\end{eqnarray} 
where $\|{\cal H}(x)\|:=\sum\limits_{k=0}^n |{\cal H}_j(x)\|$, the 
${\cal H}_j$ being defined as the ${\cal X}$-homogenizations of the $H_j$, 
namely 
$$
{\cal H}_j(x_1,\dots,x_s):= 
\Bigg(\prod\limits_{i=1}^s x_i^{-\min\limits_{\xi\in \Delta} <\xi,\eta_i>}
\Bigg) H_j(x^{<e_1^*,\eta>},\dots,x^{<e_n^*,\eta>})\,.
$$
The fact that $\delta_1+\cdots+\delta_n$ is $n$-dimensional and that 
$\beta$ lies in the interior of this polyhedron implies that one has 
for any $i=1,\dots,s$, 
$$
<\beta,\eta_i>-1-\min\limits_{\xi \in \delta_1+\dots+\delta_n}
<\xi,\eta_i>=-\min\limits_{\xi \in \delta_1+\dots+\delta_n}
<\xi-\beta,\eta_i>-1 
\geq 0\,.
$$
Therefore, one has, for any $i=1,\dots,s$, 
\begin{eqnarray}
<\beta,\eta_i>-1+\sum\limits_{j=1}^n \nu_{ij}&=&
<\beta,\eta_i>-1-\sum\limits_{j=1}^n \min\limits_{\xi \in 
\Delta+\delta_j} <\xi,\eta_i>\nonumber \\
&=&<\beta,\eta_i>-1-\sum\limits_{j=1}^n \min\limits_{{\xi \in 
\Delta}\atop {\xi_j\in \delta_j}} (<\xi+\xi_j,\eta_i> \nonumber \\ 
&\geq & -n\min\limits_{\xi\in \Delta} <\xi,\eta_i> \ \geq 0
\end{eqnarray}
since $\Delta$ contains the origin; note that any number $\nu_{ij}-\mu_{ij}$
($i=1,...,s$, $j=1,...,n$) 
is also nonnegative since $\Delta+\delta_j$ contains $\Delta_j$ for any 
$j=1,\dots,s$.  
\vskip 1mm
\noindent
Before going on in the proof of our theorem, let us point out here 
a geometric interpretation of our properness condition (1.9), in the 
spirit of \cite{elaz:gnus}. 
Let $\widetilde {\cal D}_j$, $j=1,\dots,n$, be the Cartier 
divisors on ${\cal X}$ defined as 
$$
\widetilde {\cal D}_j={\rm div} (F_j)+E(\Delta+\delta_j)\,,
$$
where $E(\Delta+\delta_j)$, $j=1,\dots,n$, is the ${\bf T}$-Cartier 
divisor which corresponds to the polyhedron $\Delta+\delta_j$; since 
$\Delta_j \subset \Delta+\delta_j$, $\widetilde {\cal D}_j$ is an 
effective Cartier divisor on ${\cal X}$. We note $E(\Delta)$ the 
${\bf T}$-Cartier effective divisor on ${\cal X}$ which corresponds to 
the polyhedron $\Delta$. Let $x$ be a point in ${\cal X}\setminus 
{\bf T}^n$ which lies in the intersection of the supports of 
the divisors $\widetilde {\cal D}_j$, $j=1,\dots,n$, and 
$V_x$ an arbitrary small neigborhood of $x$ in ${\cal X}$. Let 
${\cal I}_{V_x} \subset {\cal O}_{V_x}$ be the ideal sheaf in ${\cal O}_{V_x}$ 
which is generated by $\tilde f_{x1},\dots, \tilde f_{xn}$, where 
$\tilde f_{xj}$ is a global section in $V_x$ for the effective Cartier divisor 
$\widetilde {\cal D}_j$. Let $[E_x]$ be the exceptional divisor in the 
normalized blow-up $\pi:\ {\cal N}_x\mapsto V_x$ of $V_x$ along 
${\cal I}_{V_x}$, 
$$
[E_x]=\sum\limits_l r_{xl} E_{xl}\,,
$$
the $E_{xl}$ being its irreducible components and the $r_{xl}$ the associated 
multiplicities. Let $f_x$ be a global section for $E(\Delta)$ in $V_x$. 
As it can be seen from the properness condition (1.9) (rewritten as 
(3.18) in ${\bf T}^n$ and extended to 
the neighborhood $V_x$, which is possible since ${\bf T}^n$ is dense 
in ${\cal X}$), the order of vanishing $\rho_{xl}$ of $\pi^* f_x$ along any 
$E_{xl}$ is such that $\rho_{xl}\geq r_{xl}$ for any $l$, that is, if 
$[\pi^* f_x]$ denotes the Cartier divisor associated to $\pi^* f_x$, 
\begin{eqnarray}
[\pi^* f_x]\geq [E_x]\,. 
\end{eqnarray}
This provides the geometric interpretation we announced 
for the properness condition (1.9). 
\vskip 1mm
\noindent
This geometric vision of our properness condition being settled, 
the proof of Theorem 1.2 follows exactly the same lines than the proof 
of our previous result Theorem 1.1. We consider the complement of $U$ 
in ${\cal X}$ as 
a $2n$-chain $\Sigma$ in ${\cal X}$, and we deduce from (3.10) that
\begin{eqnarray}
&{\rm Res}\left[
\begin{array}{ccccc}
X_1^{\beta_1}\cdots X_n^{\beta_n} dX\\
F_1,\dots ,F_n \\
\end{array}
\right]_{\bf T}=-\gamma_n 
\left[\int_{\partial \Sigma} \Theta_{F,\beta,\lambda}^{\delta,\Delta}
\right]_{\lambda=0}
=-\gamma_n\left[\int_{\Sigma} \overline\partial 
\Big(\Theta_{F, \beta,\lambda}^{\delta,\Delta}\Big)
\right]_{\lambda=0}\,,\nonumber\\
&
\end{eqnarray}  
the notation $[\ ]_{\lambda=0}$ meaning that one takes the 
meromorphic continuation, and later on, the value at $\lambda=0$. 
In order to prove the vanishing of the residue symbol, it is enough to 
show that if $x$ is any point in ${\cal X}\setminus {\bf T}^n$, $V_x$ 
an arbitrary small neighborhood of $x$ in ${\cal X}$, and $\varphi 
\in {\cal D} (V_x)$, then the function 
\begin{eqnarray}
\lambda \mapsto \int_{V_x} \varphi \overline\partial\Big( 
\Theta_{F,\beta,\lambda}^{\delta,\Delta}\Big)
\end{eqnarray}
can be continued as a meromorphic function of $\lambda$ which has 
no pole at $\lambda=0$ and vanishes at $\lambda=0$. In order to do that, we 
repeat the argument in the proof of Theorem 1.1 and use a resolution 
of singularies ${\cal Y}\buildrel {\pi}\over{\rightarrow} 
V_x$, followed be toroidal resolutions 
$\widetilde {\cal Y}_\omega \buildrel {\tilde 
\pi_\omega}\over{\rightarrow} 
\omega$ over each local chart $\omega$ 
on ${\cal Y}$, such that in local coordinates $(t_1,\dots,t_n)$ on 
a local chart $\varpi$ in some $\widetilde {\cal Y}_\omega$, one has 
$$
\tilde \pi^*_\omega\pi^*  \tilde 
f_{xj}=\tilde u_j(t)t_1^{\tilde \theta_{j1}}
\cdots t_n^{\tilde \theta_{jn}}=\tilde u_j(t)\tilde m_j(t),\ 
$$
and some of the $\tilde m_j$, let say $\tilde m$, divides 
$\tilde m_1,\dots,\tilde m_n$. For the same reasons that lead to (3.20)
from the properness condition (1.9) when one was using a normalized 
blow up instead of the tower of resolutions 
$\widetilde {\cal Y}_\omega \buildrel {\tilde 
\pi_\omega}\over{\rightarrow} 
\omega \buildrel {\pi}\over{\rightarrow} 
V_x$, one can see that the properness condition implies 
that $\tilde m$ divides 
$\tilde \pi^*_\omega\pi^* f_x$ 
in the local chart $\varpi$. Therefore, it follows from 
(3.13) and (3.19) that for any test function $\tilde \varphi$ 
on ${\cal Y}$ with support in $\omega$, one can write in $\varpi$, 
$$
\tilde\pi^*_\omega [\tilde 
\varphi [\pi^*\big(\varphi \overline\partial 
\Theta_{F,\beta,\lambda}^{\delta,M}\big)]](t)
=\lambda |\tilde m|^{2\lambda} |\xi|^{2\lambda}
 \Bigg({\overline {\partial \tilde m}\over 
\overline {\tilde m}}\wedge \sigma(t) +\tau(t)\Bigg)
$$
where $\sigma$ and $\tau$ (depending on $\varphi$ and $\tilde \varphi$) 
are smooth differential forms with respective type $(n,n-1)$ and $(n,n)$ and 
$\xi$ is a real analytic strictly positive function in $\varpi$. Therefore, 
$$
\lambda\mapsto 
\tilde\pi^*_\omega [\tilde 
\varphi [\pi^*\big(\varphi \overline\partial 
\Theta_{F,\beta,\lambda}^{\delta,M}\big)]](t) 
$$ 
can be continued as a distribution valued meromorphic map on $\varpi$, which 
has no pole at $\lambda=0$ and vanishes at this point. Since the 
meromorphic function (3.22) is expressed as a sum of functions of the 
form 
$$
\lambda \int_\varpi \tilde \psi \tilde\pi^*_\omega [\tilde 
\varphi [\pi^*\big(\overline\partial 
\Theta_{F,\beta,\lambda}^{\delta,M}\big)]](t)
\ [\pi^*_\omega \pi \varphi](t)\,,  
$$
the vanishing of the residue symbol (3.21) follows. This concludes the 
proof of our theorem 1.2. $\quad\diamondsuit$ 
\section{ 
Some applications of Vanishing  Theorems for global sums of Residues.}
\setcounter{equation}{0}       
The generalized Jacobi Theorems 1.1, 1.2 derived above  
have as a direct consequence the following nonstandard 
formulations of Cayley-Bacharach 
type theorems in the spirit of \cite{gh:gnus}.
\begin{theorem}
Let ${\cal X}_j=\{{\cal P}_j=0\}$, $j=1,\dots n $ be $n$ 
hypersurfaces in ${\bf P}^n $ defined by  
homogeneous polynomials ${\cal P}_j$
in $n+1$ variables. Let $P_j(X_1,\dots,X_n):
={\cal P}_j(1,X_1,\dots,X_n)$, $j=1,\dots,n$, and 
assume that the mapping $(P_1, \dots ,P_n)$ is such that there exist 
constants $c>0, R>0$, and rational numbers 
$0<\delta_i \leq \deg(P_i)$, $i=1,\dots , n$, so that  
the properness condition $(1.7)$ holds. Suppose 
also that the common zeroes of $P_1,\dots,P_n$ in 
${\bf C}^n$ are all simple. Let 
${\cal Z}$ be the 
affine algebraic variety 
${\cal Z}:= {\cal X}_1\cap \cdots \cap {\cal X}_n\setminus \{X_0=0\}$ in 
${\bf C}^n$ 
and ${\cal Y}$ any hypersurface in ${\bf P}^n$ with degree 
less or equal than $\delta_1+\cdots+\delta_n-n-1$. Then it is 
impossible that ${\cal Y}$ contains all points of ${\cal Z}$ but one 
without containing all of them. 
\end{theorem}

\noindent
{\bf Proof.} Suppose ${\cal Y}=\{{\cal Q}=0\}$ and let 
$Q(X_1,\dots,X_n)={\cal Q}(1,X_1,\dots,X_n)$. 
Recall that 
\begin{eqnarray*}
{\rm Res}_{\alpha}\, \left [
\begin{array}{ccccc}
Q(\zeta_1,\cdots, \zeta_n)d\zeta\\
P_1,\cdots ,P_n\\
\end{array}
\right ]={Q(\alpha )\over {\cal J}_P(\alpha )},
\end{eqnarray*}
where ${\cal J}_P(\alpha )$ is the value of the Jacobian of the $(\delta _1\dots 
,\delta _n )$-proper mapping $P$ at
the simple common zero $\alpha$. The hypothesis on the degree 
of the hypersurface implies (if one uses theorem 1.1) that 
\begin{eqnarray*}
{\rm Res}\, \left [
\begin{array}{ccccc}
Q(X_1,\cdots, X_n)dX\\
P_1,\cdots ,P_n\\
\end{array}
\right ] 
=0.
\end{eqnarray*}
Therefore, if $Q$ vanishes at all points in ${\cal Z}$ but one, it 
vanishes in fact at any point in ${\cal Z}$. $\quad\diamondsuit$ 
\vskip 2mm
\noindent
We may also state a toric version of a Cayley-Bacharach theorem.
\begin{theorem}
Let ${\cal X}_j,\ j=1,\dots n$, be $n$  hypersurfaces in ${\bf C}^n$, 
defined by sparse algebraic equations $F_j(\zeta)=0$, $j=1,\dots,n$. Let 
$\Delta_j$ 
the Newton polyhedron of the polynomial $F_j$ (considered as 
a Laurent polynomial). Suppose that there exist convex polyhedra 
$\delta_1,\dots, \delta_n$, with vertices in ${\bf Q}^n$, such that
$\delta_j\subset \Delta_j$, ${\rm dim}\, (\delta_1+\cdots+\delta_n)=n$, 
and the condition (1.9) holds. Suppose also that $F_1,\dots,F_n$ 
define only simple common zeroes in $({\bf C}^*)^n$. Then, any 
hypersurface ${\cal Y}$ in ${\bf C}^n$  
which is defined as $\{Q=0\}$, where the Newton
polyhedron of $Q$ lies in the interior of $\delta_1+\cdots+\delta_n$, and 
contains any point in ${\cal Z}:=
{\bf T}^n\cap \{F_1=\dots=F_n=0\}$ but one, contains 
in fact all points in ${\cal Z}$.     
\end{theorem}

\noindent
{\bf Proof.} 
The proof is a direct application of the Theorem 1.2, exactly as our previous 
result follows from Theorem 1.1.
$\quad\diamondsuit$
\vskip 2mm
\noindent
Finally we can state an application of theorem 
1.1 (resp. 1.2) to some effective version of division problems with respect 
to proper quasi-regular maps. In the 
first case, this version is the key ingredient for a general 
explicit formulation to the algebraic Nullstellensatz \cite{by1:gnus, 
by2:gnus}; 
we do not know yet  
if the same holds in the toric case for the Newton Nullstellensatz. 
\vskip 2mm
\begin{prop} 
Let $P:=(P_1,...,P_n)$ be a $(\delta_1,\dots,\delta_n)$- proper polynomial map 
from ${\bf C}^n$ to ${\bf C}^n$, where 
$\delta_j>0$ for any $j$; suppose that $\deg P_j=D_j$, $j=1,\dots,n$.  
Let $Q_{jk}$, $j, k=1,\dots,n$ be polynomials in $(X_1,\dots,X_n, Y_1,
\dots,Y_n)$ such that $\deg Q_{jk}
\leq D_j-1$, $j=1,\dots,n$, and 
$$
P_j(Y)-P_j(X)=\sum\limits_{k=1}^n Q_{jk}(X,Y) (Y_k-X_k)\,.
$$
Let 
$$
\det \Big [Q_{jk}(X,Y)\Big ]_{1\leq j,k\leq n}
=\sum\limits_{{\alpha,\beta \in {\bf N}^n}
\atop {|\alpha|+|\beta|\leq D_1+\cdots+D_n-n}}
\gamma_{\alpha,\beta} X^{\alpha} Y^{\beta}\,.
$$ 
Then for any polynomial $Q$ with degree $D$, one has 
the following identity
\begin{eqnarray}
& Q(Y)= \nonumber \\
&=\sum\limits_{{\alpha,\beta \in {\bf N}^n}
\atop {|\alpha|+|\beta|\leq D_1+\cdots+D_n-n}}
\; \sum\limits_{{\mu \in {\bf N}^n}
\atop {<\mu+\underline 1,\delta>\leq |\alpha|+D+n}}
\gamma_{\alpha,\beta}\; {\rm Res} \left[\matrix{
Q(X)X^\alpha dX\cr
P_1^{\mu_1+1},...,P_n^{\mu_n+1}}\right] 
Y^\beta P(Y)^\mu\,, \nonumber\\
& 
\end{eqnarray}
where we used the standard notations: $\zeta^m=\zeta_1^{m_1}\dots 
\zeta_n^{m_n}$ for $\zeta\in {\bf C}^n$ and $m\in {\bf N}^n$, 
$<m_1,m_2>=m_{11}m_{21}+\dots+m_{1n} m_{2n}$ for $m_1,m_2\in {\bf N}^n$, 
$\underline 1=(1,\dots,1)$ ($n$ times).    
\end{prop}
\vskip 2mm
\noindent
{\bf Proof.} The proof follows from the Cauchy-Weil integral 
representation formula, exactly as in \cite {by1:gnus}; the 
analytic expansion of the Cauchy kernel 
that appears in this formula 
truncates thanks to Theorem 1.1. $\quad\diamondsuit$ 
\begin{corollary}
Let $P:=(P_1,...,P_n)$ be a $(\delta_1,\dots,\delta_n)$- proper polynomial map
from ${\bf C}^n$ to ${\bf C}^n$, where 
$\delta_j>0$ for any $j$; let $Q$ be in the ideal $I(P_1,\dots,P_n)$; then, one can write 
a division formula for $Q$ respect to the ideal 
$I(P_1,\dots,P_n)$ as
\begin{eqnarray}
Q(Y)=\sum\limits_{{\mu\in ({\bf N}^n)^*,\ \nu \in {\bf N}^n}\atop 
{{|\nu|\leq D_1+\cdots+D_n-n} \atop 
{<\mu+\underline {1},\delta>+|\nu| \leq
D_1+\dots+D_n+{\rm deg} Q}}}
\tilde \gamma_{\mu,\nu} Y^\nu P(Y)^\mu 
\end{eqnarray}
\end{corollary}
\noindent
Note that if $\delta_j=D_j$ (that is the $P_j$ do not have common zeroes 
at infinity), formula (4.2) becomes 
$$
Q(Y)=\sum\limits_{{\mu\in ({\bf N}^n)^*,\ \nu \in {\bf N}^n}\atop
{<\mu,D>+|\nu| \leq
{\rm deg} Q}}
\tilde \gamma_{\mu,\nu} Y^\nu P(Y)^\mu\,,
$$
which is not a surprise since the homogenization 
${\cal Q}$ of $Q$ lies (in this case) 
in the homogeneous ideal generated by ${\cal P}_1,
\dots,{\cal P}_n$.
\vskip 2mm
\noindent
In the toric case, we need first a definition, that we recall from 
\cite {by4:gnus}, p. 454.  
\begin{definition} Let $\Delta$ be a closed convex polyhedron in ${\bf R}^n$; 
$\Delta$ is called a good polyhedron if and only if 
$$
\forall x\in \Delta\,,\  \{y\in {\bf R}^n;\, 
|y_k|\leq |x_k|,\ x_ky_k\geq 0,\ k=1,\dots,n\}
\subset \Delta\,.
$$
\end{definition}
\noindent
We can now state the toric pendant of Proposition 4.1.    
\begin{prop} Let $\delta_1$,\dots, $\delta_n$ be $n$ convex 
rational polyhedra in ${\bf R}^n$ with dimension $n$ which 
contain the origin as an interior point; let $F:= (F_1,\dots,F_n)$ 
be a system of Laurent polynomials with good Newton polyhedra 
$\Delta_1,\dots,\Delta_n$, such that $\delta_j\subset \Delta_j$ for 
any $j$ and $F$ is $(\delta_1,\dots,\delta_n)$-proper. Then one 
can find Laurent polynomials $G_{jk}$, $j,k=1,\dots,n$, in 
$(X_1,\dots,X_n,Y_1,\dots,Y_n)$, such that 
$$
\det \Big[G_{jk} (X,Y)\Big]_{1\leq j,k\leq n}
=\sum\limits_{{\alpha,\beta \in {\bf Z}^n}\atop 
{\alpha+\beta \in \Delta_1+\cdots+\Delta_n}}
\gamma_{\alpha,\beta} X^\alpha Y^\beta\,
$$
and 
$$
F_j(Y)-F_j(X)=\sum\limits_{k=1}^n G_{jk} (X,Y) (Y_k-X_k),\ 
X,\, Y\in {\bf T}^n,\ j=1,\dots,n\,.
$$   
Moreover, for any Laurent polynomial $G$ with convex polyhedron 
$\Delta$, one has 
the following algebraic identity
\begin{eqnarray}
& G(Y)= \nonumber \\
&\sum\limits_{{\alpha,\beta \in {\bf Z}^n\cap (\Delta_1+\cdots+\Delta_n)}
\atop {\alpha+\beta \in \Delta_1+\cdots+\Delta_n}}
\, \sum\limits_{{\mu \in {\bf N}^n}
\atop {\Delta+ \alpha+\underline 1\not\subset{\rm int} 
(<\mu+\underline 1,\delta>)}}
\gamma_{\alpha,\beta}\, {\rm Res} \left[\matrix{
G(X)X^{\alpha+\underline 1} dX\cr
F_1^{\mu_1+1},...,F_n^{\mu_n+1}}\right]_{\bf T} 
Y^\beta F(Y)^\mu\,,\nonumber\\
&
\end{eqnarray}
where $<m,\delta>:=m_1 \delta_1+\dots+m_n\delta_n$ for any $m\in {\bf N}^n$. 
\end{prop}    
\noindent
{\bf Proof.} For the 
construction of the $G_{jk}$ 
under the hypothesis that all $\Delta_j$ are 
good, we refer to \cite{by4:gnus}. 
The fact that one can get the 
algebraic identity (4.3) is based on the use of Cauchy-Weil formula, 
as in the proof of Proposition 4.1; for more 
details see \cite{yg2:gnus}, section 2. The development of the 
Cauchy kernel as a geometric progression truncates 
(as claimed in (4.3)) 
if one applies Theorem 1.2. $\quad\diamondsuit$
\begin{corollary} 
Let $(F_1,\dots,F_n)$ be a $(\delta_1,\dots,\delta_n)$-proper 
system of Laurent polynomials; suppose that all $\delta_j$ are $n$ 
dimensional and contain the origin as an interior point; denote as $\Delta_j$ 
the smallest good polyhedron containing the support of $F_j$, $j=1,\dots,n$. 
Then, whenever $G$ is a Laurent polynomial with Newton polyhedron 
$\Delta$ that lies in the ideal generated by $F_1,\dots,F_n$ in 
${\bf C}[X_1,\dots,X_n,X_1^{-1},\dots,X_n^{-1}]$, one can write a division 
formula for $G$ respect to $(F_1,\dots,F_n)$ as 
$$
G(Y)=
\sum\limits_{{\mu\in ({\bf N}^n)^*,\ \nu \in {\bf Z}^n\cap 
(\Delta_1+\cdots+\Delta_n) }\atop
{\Delta+\Delta_1+\cdots+\Delta_n\ \not\subset
\ {\rm int}\; <\mu+\underline 1,\delta>+\beta-\underline 1}}
\tilde \gamma_{\mu,\nu} Y^\nu F(Y)^\mu\,.
$$
\end{corollary}


\begin{thebibliography}{99}
\bibitem {at:gnus} M. F. Atiyah, Resolution of singularities and 
division of distributions, Comm. Pure. Appl. Math. 23 (1970), 145-150.
\bibitem{bc:gnus} V. Batyrev, D. Cox, On the Hodge structure of projective 
hypersurfaces in toric varieties, Duke J. Math. 75 (1994), 293-338. 
\bibitem{be:gnus} D. Bernstein, The number of roots of a system of 
equations, Funct. Anal. Appl. 9 (1975), no. 2, 183-185. 
\bibitem{bs:gnus} J. Brian\c con, H. Skoda, Sur la cl\^oture int\'egrale 
d'un id\'eal de germes de fonctions holomorphes en un point de ${\bf C}^n$, 
C. R. Acad. Sci. Paris, 278 (1974), 949-951. 
\bibitem{by1:gnus} C. A. Berenstein, A. Yger, 
{\rm A formula of Jacobi and its consequences}, Annales de l'\'Ecole Normale 
Sup\'erieure, 24, 1991, 73-83.  
\bibitem{by2:gnus} C. A. Berenstein, A. Yger, 
{Residue Calculus and effective Nullstellensatz}, to appear in the 
American Journal of Mathematics. 
\bibitem{by3:gnus} C. A. Berenstein, A. Yger, 
{\rm Residues and  effective Nullstellensatz}. ERA of the AMS, 2,2, 1996.
\bibitem{by4:gnus} C. A. Berenstein, A. Yger, Multidimensional 
residues and complexity problems, Mathematics and Computers 
in Simulation 42 (1996), 449-457. 
\bibitem{ccd:gnus} E. Cattani, D. Cox, A. Dickenstein,
{\rm Residues in toric varieties}. Compositio Math. 108 (1997), p. 35-76.
\bibitem{cd:gnus} E. Cattani, A. Dickenstein,
{\rm A global view of residues in the torus}, Journal of Pure and 
Applied Algebra 117\& 118 (1997), p. 119-144.
\bibitem{cox1:gnus} D. Cox, 
{\rm The homogeneous coordinate ring of a toric variety}, Journal of 
Algebraic Geometry. 4(1995),p. 17-50.
\bibitem {cox2:gnus} D. Cox, {\rm Recent developments in toric 
geometry},  Proc. Sympos. Pure Math. 62. Part 2. 
Amer. Math. Soc., Providence, RI, 1997.  
\bibitem {elaz:gnus} L. Ein, R. Lazarsfeld, 
  {\rm A Geometric Effective Nullstellensatz}, preprint, {\bf math.ag}/9810004.
\bibitem{egh:gnus} D. Eisenbud, M. Green, J. Harris, {\rm Cayley-Bacharach 
theorems and conjectures},  Bull. Amer. Math. Soc. 
(N.S.) 33 (1996), no. 3, 295-324.   
\bibitem {fpy:gnus} A. Fabiano, G. Pucci, A. Yger, 
 {\rm Effective Nullstellensatz and geometric 
degree for zero-dimensional ideals}, Acta Arithmetica, 78(2), 1996, 165-187. 
\bibitem{fu:gnus} W. Fulton, {\rm Introduction to toric varieties}, 
Princeton Univ. Press, 1993.
\bibitem{gh:gnus} P. Griffiths, J. Harris, {\rm Principles 
of algebraic geometry}, Wiley Interscience, New York, 1978. 
\bibitem{j:gnus} C. G. J. Jacobi, 
De relationibus, quae locum 
habere denent inter puncta intersectionis 
duarum curvarum vel trium superficierum algebraicarum 
dati ordinis, simul cum enodatione paradoxi algebraici, 
Gesammelte Werke, Band III, 329-354.   
\bibitem{ka1:gnus} B. Ja. Kazarnovskii, 
{\rm On the zeroes of exponential sums }, Soviet. math. Doklady 23 (1981), 
347-351.
\bibitem{ka2:gnus} B. Ja. Kazarnovskii, 
{\rm Newton polyhedra and zeros of systems of exponential sums}, 
Funct. Anal. Appl. 18(1984), no. 4, 299-307.  
\bibitem{kh:gnus} A. Khovanskii, 
{\rm Newton polyhedra and the Euler-Jacobi formula}, 
Russian Mathematical Surveys 33 (1978), 237-238.
\bibitem{kh1:gnus} A. Khovanskii, 
{\rm Newton polyhedra and the toroidal varieties}, 
Funct. Anal. Appl. (1978), 289-295.
\bibitem{ko1:gnus} J. Koll\'ar, 
{\rm Sharp effective Nullstellensatz}, JAMS, 1, 1988, 963-975.
\bibitem {ko2:gnus} J. Koll\'ar, {\rm Effective Nullstellensatz}, preprint, 
{\bf math.ag}/9805091.
\bibitem {lip:gnus} J. Lipman, {\rm Residues and traces of differential forms 
via Hoschschild homology}, 
Contemporary Mathematics 61, American Mathematical Society, Providence, 1987.
\bibitem {pel:gnus} T. Pell\'e, {\rm Identit\'es de B\'ezout pour certains 
syst\`emes de sommes d'exponentielles}, Ark. Math. 36 (1998), 
131-162.    
\bibitem{pty:gnus} M. Passare, A. Tsikh, A. Yger,  
{\rm Residue currents of the Bochner-Martinelli type}, 
preprint, January 1998, Univ. of Bordeaux series.  
\bibitem{re:gnus}  A. Tsikh, {\rm Multidimensional Residues and Their 
Applications}, Transl. of AMS 103, 1992.
\bibitem{va:gnus} A. Varchenko, {\rm Newton Polyhedra and estimation of 
oscillating integrals}, Funct. Anal. Appl. 10(1976), p.175-196.
\bibitem{y:gnus} A. Yger, {\rm Lectures at Croce di Magara, June 1998}, notes. 
\bibitem{yg2:gnus} A. Yger, {\rm R\'esidus, courants r\'esiduels et courants de 
Green}, in G\'eom\'etrie complexe, F. Norguet, S. Ofman, J. J. 
Szczeciniarz ed., Actualit\'es Scientifiques et Industrielles 1438, 
Hermann, 1996. 
\end{thebibliography}
\end{document}